\documentclass[11pt]{amsart}

\usepackage{mathrsfs,mathtools}
\usepackage[all]{xy}
 \usepackage{mathbbol}

\usepackage{amsthm,amsmath, amssymb}
\usepackage{enumitem}
\usepackage{geometry}
\usepackage{caption}
\usepackage{graphicx}
\usepackage{hyperref}
\usepackage{fontenc}
\usepackage{inputenc}
\usepackage{multicol}

\theoremstyle{plain}
\newtheorem{theorem}{Theorem}
\newtheorem{lemma}[theorem]{Lemma}
\newtheorem{corollary}[theorem]{Corollary}
\newtheorem{proposition}[theorem]{Proposition}
\theoremstyle{definition}
\newtheorem{example}[theorem]{Example}
\newtheorem{definition}[theorem]{Definition}
\newtheorem{remark}[theorem]{Remark}
\newtheorem{question}{Question}
\newtheorem*{ack}{Acknowledgement}
\newtheorem{conjecture}[question]{Conjecture}

\newcommand\pf{\begin{proof}}
\newcommand\epf{\end{proof}}

\renewcommand{\_}[1]{_{\left( #1 \right)}}

\newcommand{\fusion}[3]{N^{#1, #2}_{#3}}
\newcommand{\xrightarrowdbl}[2][]{%
\xrightarrow[#1]{#2}\mathrel{\mkern-14mu}\rightarrow}
\newcommand{\trid}{\triangleright}
\newcommand{\varcoprod}{\mathop{\hbox{$\coprod$}}\limits}
\newcommand{\ad}{\operatorname{ad}}
\newcommand{\Alg}{\operatorname{Alg}}
\newcommand{\ag}{\mathfrak{a}}
\newcommand{\bg}{\mathfrak{b}}
\newcommand{\superqa}[3]{{\textbf A}_{#1}(#2|#3)}
\newcommand{\toba}{\mathscr{B}}
\newcommand{\wtoba}{\widetilde{\mathscr{B}}}
\newcommand{\htoba}{\widehat{\toba}}
\newcommand{\Bc}{\mathcal{B}}
\newcommand{\Cb}{\mathbb C}
\newcommand{\Cr}{\mathscr{C}}
\newcommand{\Cg}[1]{C_{\Gb}(#1)}
\newcommand{\car}{\operatorname{char}}
\newcommand{\co}{\operatorname{co}}
\newcommand{\corep}{\operatorname{corep}}
\newcommand{\D}{\mathcal{D}}

\newcommand{\End}{\operatorname{End}}
\newcommand{\Fr}{\mathscr{F}}
\newcommand{\Fc}{\mathcal{F}}
\newcommand{\Gb}{\mathbb G}
\newcommand{\Gc}{\mathcal G}
\newcommand{\GK}{\operatorname{GK-dim}}
\newcommand{\g}{\mathfrak g}
\newcommand{\gl}{\mathfrak{gl}}
\newcommand{\GL}{\operatorname{GL}}
\newcommand{\gr}{\operatorname{gr}}
\newcommand{\Gr}{\operatorname{Gr}}
\newcommand{\Hc}{\mathcal{H}}
\newcommand{\Htt}{\mathtt H}
\newcommand{\Hb}{\mathbb H}
\newcommand{\hdual}[1]{\mathscr{H}\left( #1 \right)}
\newcommand{\Hom}{\operatorname{Hom}}
\newcommand{\Hopf}{\texttt{Hopf}}
\newcommand{\I}{\mathbb I}
\newcommand{\Ic}{\mathcal{I}}
\newcommand{\id}{\operatorname{id}}
\newcommand{\Ind}{\operatorname{Ind}}

\newcommand{\Irr}{\operatorname{Irr}}
\newcommand{\irr}{\operatorname{irr}}
\newcommand{\Jac}{\operatorname{Jac}} 
\newcommand{\K}{\mathcal{K}}
\newcommand{\rsd}{\operatorname{res}}
\newcommand{\ku}{\Bbbk}
\newcommand{\kut}{\ku^{\times}}
\newcommand{\Lie}{\operatorname{Lie}}
\newcommand{\Lc}{\mathcal L}
\newcommand{\lgo}{\mathfrak  l}
\newcommand{\Mg}{\mathfrak M}
\newcommand{\cM}{\mathcal M}
\newcommand{\Mt}{\mathtt M}
\newcommand{\lcomod}[1]{{}^{#1}\hspace{-1pt}{\mathcal M}}
\newcommand{\lmod}[1]{{}_{#1}\hspace{-1pt}{\mathcal M}}
\newcommand{\ngo}{\mathfrak n}
\newcommand{\Ng}{\mathfrak N}
\newcommand{\N}{\mathbb N}
\newcommand{\Oc}{\mathscr O}
\newcommand{\Orb}{\mathcal O}
\newcommand{\Prim}{\mathcal P}
\newcommand{\pre}{\mathfrak{Pre}}
\newcommand{\pref}{\mathfrak{Pre}_\textrm{fGK}}
\newcommand{\post}{\mathfrak{Post}}
\newcommand{\postf}{\mathfrak{Post}_{\textrm{fGK}}}
\newcommand{\postfh}[1]{\mathfrak{Post}_{\textrm{fGK}, #1}}
\newcommand{\Q}{\mathbb Q}

\newcommand{\bqg}{\mathfrak{q}}

\newcommand{\QG}{\texttt{QG}}
\newcommand{\Rc}{\mathcal R}
\newcommand{\rcomod}[1]{{\mathcal M}^{#1}}
\newcommand{\rep}{\operatorname{rep}}
\newcommand{\rk}{\operatorname{rk}}
\newcommand{\Ss}{\mathcal S}
\newcommand\soc{\operatorname{soc}}

\newcommand{\supp}{\operatorname{supp}}
\newcommand{\Tb}{\mathbb T}
\newcommand{\Uc}{\mathcal{U}}
\newcommand{\Vc}{\mathcal{V}}
\newcommand{\Wc}{\mathcal{W}}
\newcommand{\yd}[1]{{}^{ #1 }_{ #1 }\mathcal{YD}}
\newcommand{\Z}{\mathbb Z}

\begin{document}

\title[On infinite-dimensional Hopf algebras]{On infinite-dimensional Hopf algebras}

\author[Nicol\'as~Andruskiewitsch]{Nicol\'as~Andruskiewitsch}

\address[N.Andruskiewitsch]{Facultad de Matem\'atica, Astronom\'ia y F\'isica,
	Universidad Nacional de C\'ordoba. CIEM -- CONICET. \newline
	Medina Allende s/n (5000) Ciudad Universitaria, C\'ordoba, Argentina}
\email{nicolas.andruskiewitsch@unc.edu.ar}

\thanks{
\newline This work was partially supported by~CONICET {\small (PIP 11220200102916CO)},
FONCyT-ANPCyT {\small (PICT-2019-03660)} and Secyt (UNC). 
I am  grateful to the Max Planck Institute for Mathematics in Bonn for its hospitality and financial support
during my visit in April and May 2023.}


\dedicatory{Dedicated to Pierre Cartier}

\begin{abstract}
This is a survey on pointed Hopf algebras with finite Gelfand-Kirillov dimension and  related aspects
of the theory of infinite-dimensional Hopf algebras.
\end{abstract}

\maketitle

\section*{Introduction}\label{section:introduction}
According to Drinfeld the category of quantum groups is the opposite to the category of Hopf algebras 
with bijective antipode \cite[p. 800]{Drinfeld}.
What should be the right formulation of the notion of affine algebraic quantum group? 
The concept of affine quantum variety is far from being understood,
and we can approach the question from different angles.
First, we can impose finiteness assumptions as our working hypothesis; 
in short,  we mean studying Hopf algebras that are either
Noetherian  or have finite  Gelfand-Kirillov dimension ($\GK$ for short). 
It is also natural to narrow down to affine Hopf algebras, i.e., finitely generated as algebras.
The two competing notions mentioned above have their peculiarities:  finiteness of the $\GK$ 
seems to be more restrictive than Noetherianness and perhaps by this reason it  appears to be more accessible
(but by no means easy), while the former looks more natural but even in the classical families it raises 
formidable open problems that prevent us from moving forward in the formulation of appropriate strategies. 
Second, there are conditions with topological or geometric flavor that are often imposed to the Hopf algebras in consideration.
Details of some of these properties and some relevant results are given in Section \ref{section:Properties}. More information and references can be found in the surveys \cite{Brown-survey-PI,Brown-Goodearl,Brown-Gilmartin-survey,Brown-survey,Brown-Zhang-survey,Goodearl-survey}.
Third, one may focus on Hopf  algebras with suitable structural properties allowing 
to use the methods developped by (or in collaboration with) H.-J. Schneider, I. Heckenberger, I. Angiono, J. Cuadra and others, with crucial   reduction to Nichols algebras. The main goal of this paper 
is to discuss (see Section \ref{section:pointed})  the advances on the
classification of  pointed Hopf algebras with finite Gelfand-Kirillov dimension using these methods. 
More applications and generalizations of them are given in Section \ref{section:more}
where we speculate on a few aspects; some of these speculations seem to be new.

The finite-dimensional Hopf algebras, i.e., the finite quantum groups,
are the affine Hopf algebras with $\GK = 0$.  There are a number of surveys on them
\cite{A-Bariloche,A-Schneider-cambr,A-icm,A-Leyva,A-Angiono-Diag-survey}
so that here we do not repeat details and references that can be found  loc. cit.
The paper and the reference list are not meant to be exhaustive, reflecting mainly 
the interests and the limitations of the  author. Among other topics, 
the free quantum permutation groups and their relatives
intensively studied  in recent decades are excluded from our considerations; see the book \cite{Banica}.
Nor do we look at the analogous question in the framework of multiplier Hopf algebras,  \cite{Vandaele}.
Another important notion absent here is the homological integral, 
see  \cite{Brown-Zhang-survey,Lu-Wu-Zhang}.

We collect from different sources a number of questions in Section \ref{section:questions}, of various levels:
some are well-known open problems, some are of technical nature to decide the course to continue,
some are intriguing curiosities. \footnote{My favorite is Question \ref{question:CQG-Lie-exhaustion}.}

\section{Preliminaries}\label{section:Preliminaries}

\subsection{Notations}

The set of natural  numbers is denoted by $\N$, and $\N_0=\N\cup \{0\}$. 
Given $m< n \in\N_0$, we set $\I_{m,n}=\{i \in \N_0: m \leq i \leq n\}$, $\I_{n}=\I_{1, n}$.
We write $M\leq N$ to express that 
$M$ is a subobject of $N$ in a given category. 
Given a vector space $V$ and a subspace $U$ either of $V$ or its dual $V^*$,  $U^{\perp}$ 
stands for the annihilator of $U$, either in $V^*$ or $V$. The exterior and symmetric algebras of 
$V$ are denoted by $\Lambda(V)$ and $S(V)$ respectively.

\medbreak
Except when explicitly stated otherwise, the base field $\ku$ is algebraically closed and has $\car = 0$.
All algebras are associative and unital. 
The set of algebra maps from a $\ku$-algebra $A$ to a  $\ku$-algebra $B$ is denoted by $\Alg(A, B)$. 
Let $A$ be an  algebra; $A$ is \emph{affine} means that it is finitely generated. 
A two-sided ideal is  simply  called an ideal. The subalgebra generated by $X \subset A$ is denoted by $\ku \langle X\rangle$
while $\ku \langle (x_i)_{i\in I} \vert  (r_j)_{j\in J} \rangle$ denotes the algebra presented by generators $x_i$'s and relations $r_j$'s. We assume familiarity with the notion of graded algebra;   
a graded algebra
$R = \oplus_{n\in \N_0} R^n$ is connected if $R^0 \simeq \ku$.

The category of  left $A$-modules is denoted by $\lmod{A}$;
its full subcategory of finite-dimensional ones is denoted by $\rep A$.
For $V \in \lmod{A}$, the associated representation is denoted by $\rho_V: A \to \End V$. 
The set of isomorphism classes of simple objects in an abelian category $\cM$ is denoted by $\irr \cM$; 
by abuse of notation we set $\Irr A \coloneqq \irr\lmod{A}$,  $\irr A \coloneqq \irr \rep A$. 
Assume that $A$ is commutative. Then $\irr A \simeq \Alg(A, \ku)$.
Given $V \in \rep A$ and $\chi \in \Alg(A, \ku)$, we set
\begin{align}\label{eq:decomp-commutative}
V^{(\chi)} &\coloneqq \left\{ v\in V: \exists n\in \N\,   \big(a- \chi(a)\big)^n\cdot v = 0 \quad\forall a \in A \right\},
\\ \notag
V^\chi  &\coloneqq \{v\in V: a \cdot v = \chi(a)v,  \forall  f\in H\} = \soc V^{(\chi)}.
\end{align} 
A classical result, cf. \cite[1.2.13]{dixmier}, states that
\begin{align}\label{eq:decomp-commutative-dix}
\begin{aligned}
V &= \oplus_{\chi \in \Alg(A, \ku)} V^{(\chi)}, \\
\soc V  &= \oplus_{\chi \in \Alg(A, \ku)} V^{\chi}.
\end{aligned}
\end{align}
We set $\supp V \coloneqq  \{\chi \in \Alg(A, \ku): V^{(\chi)}\neq 0\} = \{\chi \in \Alg(A, \ku): V^{\chi}\neq 0\}$. 

\medbreak All coalgebras are coassociative and counital.
Let $C$ be a  coalgebra. The category of  left, respectively right, $C$-comodules is denoted by $\lcomod{C}$, 
respectively $\rcomod{C}$; 
the full subcategory of $\lcomod{C}$ consisting of finite-dimensional objects is denoted by $\corep C$.
The coradical of $C$ is  the sum $C_0$ of all its simple subcoalgebras.
 The \emph{wedge} of two subspaces $D$ and $E$ of $C$ is 
\begin{align*}
D \wedge E = \{c\in C: \Delta(c) \in D\otimes C + C \otimes E\}.
\end{align*}
In terms of the dual algebra $C^*$,  $D \wedge E = (D^{\perp}E^{\perp })^{\perp}$ \cite[Proposition 9.0.0]{Sweedler-book}. 

 It is useful to state  facts in terms of $C^*$; for instance $C$ is a $C^*$-bimodule via
 \begin{align} \label{eq:action-dual}
 f \rightharpoonup x  &= x\_{1}  \langle f, x\_{2}\rangle, & 
  x \leftharpoonup f  &=   \langle f, x\_{1}\rangle x\_{2}, & 
 x  &\in C, f\in C^*.
 \end{align}

\subsubsection*{Hopf algebras}
See \cite{Montgomery-book,Radford-book} for general expositions.
All Hopf algebras  are supposed to have bijective antipode except in \S \ref{subsec:Antipode}. 
The notation for Hopf algebras is standard: $\Delta$ is the comultiplication,
$\varepsilon$ is the counit, $\Ss$ is the antipode; by abuse of notation the transpose of the latter 
is also denoted by $\Ss$. 
For the comultiplication and the coactions we use a variant of the Heynemann-Sweedler notation, omitting the summation symbol.
Given a Hopf algebra $H$, we  set 
\begin{align*}
G(H) &= \{x\in H -0: \Delta(x) = x \otimes x\}, & &\text{the group of grouplike elements, }
\\
\Prim(H) &= \{x\in H: \Delta(x) = x \otimes 1 +1 \otimes x\}, & &\text{the Lie algebra of primitive elements,}
\\
\Prim_{g,h}(H) &= \{x\in H: \Delta(x) = x \otimes h +g \otimes x\}, & &\text{the space  of $(g,h)$-skew primitives,}
\end{align*}
 where $g,h \in G(H)$. We shall use freely the theory of tensor categories as in \cite{EGNO}.

\subsubsection*{Groups} The neutral element of a group is often denoted by $e$.
The cyclic group of finite order $d$ is denoted by $C_d$.
The group of $d$-th roots of 1 in $\ku$ is denoted $\Gb_d$ and the subset of those 
with order $d$ by $\Gb'_d$; also $\Gb_{\infty} = \bigcup_{d\ge1} \Gb_d$, 
$\Gb'_{\infty} =\Gb_{\infty} \backslash \{1\}$. The group of multiplicative characters of
a group $\varGamma$ is denoted by $\widehat{\varGamma}$.

Given two classes of groups $\mathfrak{C}$ and $\mathfrak{D}$,
a group $\varGamma$ is said to be (in the class)  $\mathfrak{C}$-by-$\mathfrak{D}$
if it has  a normal  subgroup $N$ in $\mathfrak{C}$ such that $\varGamma / N$ belongs to $\mathfrak{D}$. 

\smallbreak
A solvable group is \emph{polycyclic} if every subgroup is finitely generated; equivalently, 
it admits a subnormal series with cyclic factors. 

\smallbreak
 A group is \emph{linear} if it is  embedable into $GL(n, \ku)$ for some $n$, or equivalently into an affine algebraic group.

\smallbreak Let $\varGamma$ be  a group.
The conjugacy class of $\gamma \in \varGamma$ is denoted by $\Orb_\gamma$ and its centralizer by $\Cg{\gamma}$. 
The FC-center $FC(\varGamma)$  of  $\varGamma$ is the union 
of all its finite conjugacy classes;
it is  a characteristic subgroup of $\varGamma$  and $Z(\varGamma) \leq FC(\varGamma)$ \cite{Passman-book}.

\subsubsection*{Simple algebraic groups}\phantomsection\label{text-def:simple-algebraic} 
Let $H$ be an affine commutative Hopf algebra; then $\Gb \coloneqq \Alg(H, \ku)$ is an affine algebraic group, so that $H \simeq \Oc(\Gb)$, the algebra of regular functions on $\Gb$.
If $V\in \lcomod{\Oc(\Gb)}$, then $V$ is a (so called rational) $\Gb$-module via 
\begin{align}\label{eq:action-rational}
\gamma\cdot v &= \langle \gamma^{-1}, v\_{-1}\rangle v\_{0},&
\gamma  &\in \Gb, \, v\in V.
\end{align}
Let $\rep \Gb$ be the full subcategory of $\lmod{\ku \Gb}$ consisting 
of rational $\Gb$-modules; we set $\irr{\Gb} \coloneqq \irr \rep \Gb$.
See the reference \cite{Jantzen} for more on algebraic groups.

\medbreak
Let $\Gb$ be a simple affine algebraic group, i.e., $\Gb$ is connected and its Lie algebra $\g$ is simple.
Fix a Cartan subalgebra of $\g$ and a choice of positive roots. 
Let $Q$ be  the root lattice and $P$ the weight lattice of $\g$.
Then the weights of the rational $\Gb$-modules form a lattice $\varLambda$, $Q \leq \varLambda \leq P$;   this gives a bijection between the set of lattices intermediate between $P$ and $Q$, and the set of isomorphism classes of simple algebraic groups with Lie algebra $\mathfrak g$. 
\phantomsection\label{text-notation:simple-algebraic} 
The extreme examples are the simply connected group $\Gb_{\text{sc}}$, with $\varLambda = P$,  and the adjoint group $\Gb_{\text{ad}}$ (which is simple as abstract group), with $\varLambda = Q$.
Then $\irr \Gb$ is parametrized by 
the cone $\varLambda^+$; if  $\lambda \in \varLambda^+$, then  the corresponding
highest weight simple module is denoted by $L(\lambda)$. We set  $d_{\lambda} = \dim L(\lambda)$
and denote the fusion rules by $\fusion{\lambda}{\mu}{\nu}$, i.e.,
\begin{align}\label{eq:simple-modules}
L(\lambda) \otimes L(\mu) &\simeq \oplus_{\nu \in \varLambda^+} L(\nu)^{\fusion{\lambda}{\mu}{\nu}},& \lambda, \mu &\in \varLambda^+.
\end{align}

\subsubsection*{Filtrations} 
We shall consider various types of filtrations, for instance:

\begin{itemize}[leftmargin=*]
\medbreak\item    Descending algebra filtrations of an algebra $A$, i.e.,
families $\Fc = \left(\Fc_n(A)\right)_{n\in \N_0}$ of subspaces  of $A$ such that $\Fc_n(A) \supset \Fc_{n+1}(A)$ for all $n\in \N_0$ and
\begin{align*}
\Fc_n(A) \Fc_m(A) &\subset\Fc_{n + m}(A), & \text{for all } n, m&\in \N_0.
\end{align*}
\end{itemize}
The graded algebra associated to the filtration is $\Gr_{\Fc} A = \oplus_{n\ge 0} \Fc_n(A) / \Fc_{n+1}(A)$. 
For accurateness, we need  to assume that $\Fc$ is Hausdorff (aka separated),  i.e.,
\begin{align}\label{eq:hausdorff}
\bigcap_{n \in \N} \Fc_n(A) &= 0.
\end{align}
For instance,  let $I$ be an ideal of  $A$ such that $\bigcap_{n \in \N} I^n = 0$. 
The filtration $\Fc^I$ by powers of $I$ is given by
$\Fc^I_0 (A)= A$, $\Fc^I_{n}(A) = I^{n}$, $n\in \N$.
We discuss  examples in \S \ref{subsec:Method-radical}.

\begin{itemize}[leftmargin=*]
\medbreak\item    Ascending coalgebra filtrations  of a coalgebra $C$, i.e.,
families $\Fr = \left(\Fr_n(C)\right)_{n\in \N_0}$ of subspaces of $C$ 
such that $\Fr_n(C) \subset \Fr_{n+1}(C)$ for  $n\in \N_0$ 
and
\begin{align*}
\Delta(\Fr_n(C) )  &\subset \sum_{0 \leq i \leq n} \Fr_i(C) \otimes \Fr_{n-i}(C), & \text{for all } n&\in \N_0.
\end{align*}
\end{itemize}
The graded coalgebra associated to the filtration is $\gr_{\Fr} C = \oplus_{n\ge 0} \Fr_n(C) / \Fr_{n-1}(C)$
where $\Fr_{-1}(C) = 0$. \footnote{We omit the subscript $\Fr$ when  no confusion is possible.}
For accurateness, we need  to assume that $\Fr$ is exhaustive,  i.e.,
\begin{align}\label{eq:exhaustive}
\bigcup_{n \in \N} \Fr_n(A) &= C.
\end{align}
 For instance,  let $D$ be a subcoalgebra of  $C$. 
 The $D$-wedge filtration $\left(\wedge^{n} D\right)_{n \in \N_0}$  is defined
recursively  by $\wedge^{0} D = D$ and $\wedge^{n+1} D = (\wedge^{n} D) \wedge D$ for $n\in \N_0$.
This filtration is exhaustive if and only if $C_0 \subseteq D$ \cite[5.3.4]{Montgomery-book}; when $C_0 = D$
this is called the coradical filtration and $\wedge^{n} C_0$ is denoted by $C_n$.
 We discuss this filtration in \S \ref{subsec:Method}.
 
\begin{itemize}[leftmargin=*]
\medbreak\item    The PBW filtration of an algebra with a PBW-basis. 
\end{itemize}
 
Let $A$ be an algebra. Let $\emptyset \neq P \subset A$ and $\emptyset \neq S \subset A$, 
$<$  a total order on $S$ and $h$ a function $h: S \mapsto \N \cup \{ \infty \}$ called the height.
If the set $B = B(P,S,<,h) =$
\begin{align*}
&= \big\{p\,s_1^{e_1}\dots s_t^{e_t}:& t &\in \N_0,\ s_i \in S, \ p \in P,& s_1&>\dots >s_t, 
&  0&<e_i<h(s_i) \big\},
\end{align*}
is a $\ku$-basis of $A$, then we say that it is a \emph{PBW-basis}  and that
the elements of $S$ are the PBW-generators. 
For simplicity, assume from now on  that $P = \{1\}$ and that $S$ is finite; fix 
a numeration $S = \{s_1,s_2,\dots, s_N\}$ such that $i<j$ iff $s_i<s_j$.
Then any $1 \neq b \in B$  can be expressed as $b = s_M^{e_M}\cdots s_1^{e_1}$ where $0 \leq e_i < h(s_i)$, $i\in \I_M$, and $e_M\neq 0$; we set
\begin{align*}
\deg b &= (e_1,\dots, e_M, 0, \dots) \in \N_0^{N},
\end{align*}
Let $\delta_{j} \in \N_0^{N}$, $\delta_j(i) = \delta_{ij}$, $j\in \I_N$. 
Let $\preceq$ be the lexicographical order, reading from the right, on  $\N_0^{(\N)}$.
This gives an exhaustive ascending $\N_0^{N}$-filtration of $A$:
\begin{align*}
A_f &\coloneqq \text{the vector subspace spanned by }  \{b \in B: \deg b \preceq f\} ,&
f  &\in \N_0^{N},
\end{align*}

The following Lemma, inspired by \cite{Deconcini-kac-procesi},  explains when it is an algebra filtration.

\begin{lemma}\label{lema:dck}  \cite{A-Angiono-Heckenberger-infinite}
The family $(A_f)_{f \in \N_0^{(\N)}}$ is an algebra filtration iff the following hold:
\begin{enumerate}[leftmargin=*,label=\textrm{(\alph*)}]
\medbreak\item   \label{item:PBW-convex-1} For every  $i <j\in \N$ there exists $\lambda_{ij} \in \ku$ such that
$s_is_j = \lambda_{ij} s_js_i + \sum_{f \prec \delta_{i} + \delta_{j}} A_f$;

\medbreak\item   \label{item:PBW-convex-2} for every  $i\in \N$ such that $h(s_i)\in\N$, 
$s_i^{h(s_i)}  \in \sum_{f \prec h(s_i) \delta_{i}} A_f$.
\end{enumerate}
\end{lemma}

When this happens, the PBW-basis $B$ is said to be \emph{convex}.
See \cite{A-Angiono-Heckenberger-infinite} for  details.

\subsection{Techniques of constructions}\label{subsec:Thecniques}
To set up the notation, we recall some standard ways of building new Hopf algebras,
referring to \cite{A-Bariloche} for details and references.
Throughout $H$ is a Hopf algebra.

\subsubsection*{Duals} If $H$ is finite-dimensional, then its
linear dual $H^{*}$ is a again a Hopf algebra, by transposing the 
structure maps, but otherwise $H^{*}$ is not a Hopf algebra.  
Instead, there is a largest Hopf algebra $H^{\circ}$ inside $H^*$ with transposed structure maps, given by
\begin{align*}
H^{\circ} &= \{f \in H^*: \text{ there exists an ideal $I$ such that } f_{\vert I} = 0  \text{ and }
\dim H/I < \infty\};
\end{align*}
this is the \emph{finite dual} of $H$. Here is an alternative description.
If $V \in \rep H$, let $C_V$ be
the image of the transpose of the associated representation 
 $^{t}\rho_V:   (\End V)^{*} \to  H^{*}$. Then $C_V$, a subcoalgebra of $H^{\circ}$,  is spanned 
by the matrix coefficients of $V$  and 
\begin{align*}
H^{\circ} &= \bigcup_{V \in \rep H} C_{V}.
\end{align*}

Let $\Cr$ be a full subcategory of $\rep H$ closed under subquotients, tensor products and duality -- that is a tensor subcategory of $\rep H$
\cite[4.11.1]{EGNO} -- but  not necessarily closed under extensions. Then the  subcoalgebra of $H^{\circ}$
\begin{align*}
\hdual{\Cr} &\coloneqq \bigcup_{W\in \Cr} C_{W}
\end{align*}
is a Hopf subalgebra of $H^{\circ}$. Notice that different Hopf subalgebras of the finite dual might behave quite differently \cite{ANT}.

\subsubsection*{Yetter-Drinfeld modules}
Recall that a  braided vector space is a pair $(V, c)$ where $V$ is a vector space
and $c\in GL(V \otimes V)$ 
satisfies the braid equation 
\begin{align*}
(c\otimes \id)(\id\otimes c)(c\otimes \id) = (\id\otimes c)(c\otimes \id)(\id\otimes c).
\end{align*} 

A  {\textit Yetter-Drinfeld module} over a Hopf algebra $H$  is a  left $H$-module 
$M$ and simultaneously a left $H$-comodule via
$\delta:M \to H\otimes M$, both structures compatible by
\begin{align*}
\delta(h.m) &= h_{(1)} m_{(-1)}\Ss(h_{(3)}) \otimes h_{(2)}. m_{(0)}, & h\in H, m \in M.
\end{align*}
The category $\yd{H}$ of   Yetter-Drinfeld modules over $H$ is a braided tensor one,   
where, for $M, N \in \yd{H}$, the braiding $c_{M,N}: M\otimes N \to N\otimes M$ is given  by
\begin{align}\label{eq:braid-yd}
c_{M,N}(m\otimes n) &= m_{(-1)}. n \otimes m_{(0)}, & m \in M,&\, n\in N.
\end{align}
Thus if $M\in \yd{H}$ then  $(M,  c_{M,M})$ is a braided vector space. 
Reciprocally a \emph{rigid} braided vector space
can be realized as Yetter-Drinfeld module (in many ways) \cite{Takeuchi}. 

\medbreak
Since  $\yd{H}$ is a braided tensor category,  Hopf algebras
in  $\yd{H}$ are available, which are  not usual Hopf algebras in general.
It is also convenient to consider also braided Hopf algebras, 
i.e., braided vector spaces with compatible multiplication, comultiplication and antipode.
Again these can be realized as Hopf algebras in categories of Yetter-Drinfeld modules in many ways, 
when they are rigid. See \cite{Takeuchi}.

\subsubsection*{Bosonization}
This is a  kind of semidirect product, particular of the theory of Hopf algebras,  
discovered by Radford and interpreted categorically by Majid. 
Given a brai\-ded Hopf algebra $R$ in $\yd{H}$, the bosonization (aka Radford biproduct)
$R\# H$ is the vector space $R\otimes H$ with the semidirect multiplication and
semidirect comultiplication; it turns out that this is a Hopf algebra. 
Besides, this construction naturally arises  as follows. 
Let $A$  be a Hopf algebra provided with Hopf algebra  maps $\pi: A\to H$ and  $\iota: H\to A$ such that $\pi\iota 
= \id_{H}$. Then the algebra of coinvariants
\begin{align*} 
R \coloneqq  A^{\co \pi} = \{a\in A: (\id\otimes \pi)\Delta (a) = a\otimes 1\}
\end{align*}
is a braided Hopf algebra in $\yd{H}$ and $A \simeq R\# H$. Some examples:

\begin{itemize}[leftmargin=*]
\medbreak\item  Let $R$ be a Hopf algebra and $\varGamma$ a group acting by Hopf algebra automorphisms on $R$. 
Provided with the trivial coaction, $R$ is a Hopf algebra in $\yd{\ku\varGamma}$ 
and $R \# \ku \varGamma$ is a Hopf algebra known as the smash product of $R$ and $H$.

\medbreak
\medbreak\item 
The symmetric category of super vector spaces is a tensor subcategory of $\yd{\ku C_2}$; 
thus any Hopf superalgebra $\Hc$ gives rise to a Hopf algebra $\Hc \# \ku C_2$.
\end{itemize}

\subsubsection*{Twisting and cocycle deformations}
The following two devices are very useful and have categorical interpretations.
See \cite{A-Bariloche} or \cite{EGNO} for more details.

\begin{itemize}[leftmargin=*]
\medbreak\item  Given  an invertible 2-cocycle $\sigma: H \otimes H \to \ku$, $H_{\sigma} = H$ 
with the new multiplication $\cdot_{\sigma}$ given by 
$x\cdot_{\sigma}y = \sigma(x\_1, y\_1) x\_2\cdot y\_2\sigma^{-1}(x\_3, y\_3)$
and the same $\Delta$ is again a Hopf algebra.

\medbreak\item  Given  an invertible dual 2-cocycle $F \in H \otimes H$ (aka twist), $H^F = H$ 
with the new comultiplication $\Delta^F \coloneqq F\Delta F^{-1}$ 
and the same multiplication is again a Hopf algebra.
\end{itemize}

\subsubsection*{Extensions}
This notion was considered by many authors, see e.g. the references in
\cite{A-Bariloche}. We say that 
a sequence of morphisms of Hopf algebras
$A\xhookrightarrow[]{\iota} C \xrightarrowdbl[]{\pi} B$
is exact if:
\begin{multicols}{2}
\begin{enumerate}[leftmargin=*,label=\textrm{(\alph*)}]
\medbreak\item \label{suc-exacta-1} $\iota$ is injective.
\medbreak\item \label{suc-exacta-2} $\pi$ is surjective.
\medbreak\item \label{suc-exacta-3} $\ker\pi = C\iota(A)^+$.
\medbreak\item \label{suc-exacta-4} $\iota(A) = C^{\operatorname{co} \pi}$.
\end{enumerate}
\end{multicols}
These requirements are simplified if $A$ is a sort of normal Hopf subalgebra 
of $C$, which is a faithfully flat $A$-module, or dually  if $B$ 
is a sort of conormal quotient Hopf algebra of $C$, which is a  faithfully coflat $B$-comodule.
See e.g. \cite{A-Canadian,Schneider-ext}.

\subsubsection*{Faithful flatness}
The previous discussion and other reasons, see e.g. \cite{Takeuchi-hopf-ideals}, 
lead to ask whether a  Hopf algebra $H$ is  faithfully flat over any Hopf subalgebra $A$
 \cite{Montgomery-book}; the answer is negative \cite{Schauenburg-non-fflat} but
 in the example, $A$ has a non-bijective antipode. \phantomsection\label{text-question:fflat}
 However it is still open when  both $H$ and $A$ have bijective antipodes,  Question \ref{question:fflat-subalg}. 
Various  partial results are known:   $H$ is  faithfully flat over a Hopf subalgebra $A$ when

\begin{itemize}[leftmargin=*]
\medbreak\item [$\circ$]    \cite[14.1]{Waterhouse} $H$ is commutative; more generally, when 
  $A$ is commutative \cite{Arkhipov-Gaitsgory};

\medbreak\item [$\circ$] \cite{Skryabin-antipode-Noetherian} $H$ is Noetherian and 
residually finite-dimensional, \footnote{An algebra is \emph{residually finite-dimensional} 
	if the intersection of its ideals of finite codimension is 0.\phantomsection\label{foot:res-fd}}
 and
$A$ is right Noetherian (without assuming bijectivity of the antipode), this leads to   Question \ref{question:fflat-subalg-Noetherian};

\medbreak\item [$\circ$] \cite{Chirvasitu} $H$ is cosemisimple.

\medbreak\item [$\circ$]  (Nichols-Z\"oller) If $\dim H< \infty$, then $H$ is  free over any Hopf subalgebra.

\medbreak\item [$\circ$] \cite{Radford-book} If $H$ is pointed, then $H$ is  free over any Hopf subalgebra.
\end{itemize}

\subsection{Classical examples}\label{subsec:Examples}
Recall that $\ku = \overline{\ku}$, $\car \ku =0$. Throughout $H$ is a Hopf algebra. By classical  we mean:

\begin{itemize}[leftmargin=*]
\medbreak\item    The enveloping algebra $U(\g)$ of a Lie algebra $\g$.

\medbreak\item    The group algebra $\ku \varGamma$ of a group $\varGamma$. 

\medbreak\item     The coordinate ring $\Oc(\Gb)$ of  an affine algebraic group $\Gb$.
\end{itemize}

These examples essentially exhaust  the Hopf algebras in two natural classes. 
The first result was obtained independently by  Cartier, Kostant and Gabriel in the 60's,
extending the celebrated Milnor-Moore Theorem.

\begin{itemize}[leftmargin=*] \medbreak\item   [$\circ$]  
If $H$ is cocommutative,
then $H \simeq U(\g)\rtimes \ku\varGamma$ where $\varGamma = G(H)$ acts on $\g = \Prim (H)$ by conjugation.
\end{itemize}

\begin{itemize}[leftmargin=*] \medbreak\item   [$\circ$] \cite{Cartier}  If $H$ is commutative,
then $H$ has trivial radical; hence 
$H \simeq \Oc(\Gb)$ where $\Gb = \Alg (H, \ku) \simeq \varprojlim \Gb_{i}$ is the limit of the
inverse system  of the algebraic groups
$\Gb_i = \Alg (H_i, \ku)$, $H_i$ running through all finitely generated Hopf subalgebras of $H$.
\end{itemize}

To this list one may add examples coming from the super world:

\begin{itemize}[leftmargin=*]
\medbreak\item     The bosonization  $\Oc(\Gb) \# \ku \, C_2$ of the coordinate ring   of a pro-affine algebraic supergroup $\Gb$.

\medbreak\item    The bosonization $U(\g) \# \ku \, C_2$ of the enveloping algebra of a Lie superalgebra $\g$. 
\end{itemize}

There is a super version of the result above due to Kostant \cite[Prop. 3.2]{Kostant}.

\begin{itemize}[leftmargin=*] 
\medbreak\item   [$\circ$]  
If $E$ is a connected cocommutative Hopf superalgebra,
then $E \simeq U(\g)$,  $\g = \Prim (E)$.
\end{itemize}

\subsection{Quantum examples}\label{subsec:Quantum-Examples}
The discovery of quantum groups revolutionized the theory of Hopf algebras.
We merely set up the notation and refer to \cite{Drinfeld,Brown-Goodearl,Deconcini-kac-procesi,Lusztig-book,Jantzen-quantum}.
Let $\Gb$ be a simple simply connected algebraic group,
$\g = \operatorname{Lie} \Gb$, $q \in \kut \backslash \Gb_{\infty}$ and $\epsilon\in \Gb_{\ell}'$, where
$\ell \in 2\N + 1$, $(\ell, 3)=1$ if $\Gb$ is of type $G_2$. There are various sorts of quantum groups:

\begin{itemize}[leftmargin=*]
\medbreak\item    
The quantized enveloping algebra $U_q(\g)$ \cite{Drinfeld,Jimbo};

\medbreak\item    
The quantized algebra  of functions $\Oc_q(\Gb)$, the Hopf subalgebra of the finite dual
${U_q(\g)}^{\circ}$ arising from the representations of type 1 \cite{Lusztig-book};

\medbreak\item  the  quantized enveloping algebra  $\mathcal U_\epsilon(\g)$ introduced by Lusztig \cite{Lusztig-jams-findim,Lusztig-dedicata};

\medbreak\item  the DCKP quantized enveloping algebra  $U_\epsilon(\g)$ introduced and studied  in \cite{Deconcini-kac-procesi,Deconcini-procesi};

\medbreak\item  the finite quantum groups $\mathfrak u_\epsilon(\g)$ \cite{Lusztig-jams-findim,Lusztig-dedicata};

\medbreak\item   the quantized algebra  of functions $\Oc_\epsilon(\Gb)$, a Hopf subalgebra of  
${\mathcal U_\epsilon(\g)}^{\circ}$  \cite{Lusztig-jams-findim}.
\end{itemize}
  
Some  of them fit into  extensions as follows, where $\Gb^d$ is Poisson dual to $\Gb$:
\begin{align*}
& \mathfrak u_\epsilon(\g) \hookrightarrow \mathcal U_\epsilon(\g) \twoheadrightarrow U_\epsilon(\g); 
&\Oc(\Gb^d) \hookrightarrow U_\epsilon(\g) \twoheadrightarrow {\mathfrak u_\epsilon(\g)}; &
&\Oc(\Gb) \hookrightarrow \Oc_\epsilon(\Gb) \twoheadrightarrow {\mathfrak u_\epsilon(\g)}^{*}.
\end{align*}
There are many ways to obtain new Hopf algebras from these examples; we cite two:
\begin{itemize}[leftmargin=*]
\medbreak\item [$\circ$]  The quantum subgroups, i.e., Hopf algebra quotients, of $\Oc_\epsilon(\Gb)$,
were classified in  \cite{ag-compo} (previous work was done by P. Podles 
and E. Müller); they are extensions
$\Oc(L) \xhookrightarrow[]{\iota} H \xrightarrowdbl[]{\pi} {\mathfrak m}^{*}$, 
where  $L\leq G$, $\mathfrak m \leq \mathfrak u_\epsilon(\g)$.

\medbreak\item [$\circ$] \cite{ANT} For any discrete (not necessarily algebraic) 
subgroup $\varGamma \leq \Gb$ there is a Hopf subalgebra
$\mathcal{H}(\varGamma)$ of $\Oc_\epsilon(\Gb)^{\circ}$ fitting in an exact sequence
$\mathfrak u_\epsilon(\g) \hookrightarrow  \mathcal{H}(\varGamma)  \twoheadrightarrow \ku \varGamma$.
\end{itemize}

\section{Properties}\label{section:Properties}
In this section we discuss some relevant structural properties and what is known or 
conjectured for the examples above. Below a filtered algebra is
an algebra $A$  with a given ascending exhaustive filtration and $\gr A$ denotes the associated 
graded algebra.

\subsection{Gelfand-Kirillov dimension}\label{subsec:GK}
The Gelfand-Kirillov dimension, $\GK$ for short,  of an associative algebra  $A$ is  
\begin{align*}
\GK A &\coloneqq \sup _{\substack{1 \in V \subset A\\ \dim V < \infty}}
\limsup _{{n\to \infty }}\log _{n}\dim _{k}V^{n}.
\end{align*}
Here $V^{n}$ is the linear span of all products of $n$ elements in $V$. See \cite{Krause-Lenagan}.
If $A$ is an affine commutative algebra, then $\GK A = \dim \operatorname{Spec} A$.
We record a useful result:

\begin{lemma}\label{lemma:filtration-GK} \cite[6.5, 6.6]{Krause-Lenagan}
Let $A$ be a filtered algebra. Then $\GK A \geq \GK \gr A$ and the equality holds if $\gr A$ is finitely generated.
\end{lemma}
As a direct consequence we note:
\begin{itemize}[leftmargin=*]
\medbreak\item   [$\circ$]  If $\g$ is a Lie algebra, then $\GK U(\g) = \dim \g \in [0, \infty]$.

\medbreak\item   [$\circ$]  If $\g$ is a semisimple finite-dimensional Lie algebra and $q$ is not a root of 1, then $\GK U_q(\g) = \dim \g$.
\end{itemize}

By a famous result of Gromov, there is also a complete answer for group algebras:

\begin{itemize}[leftmargin=*]
\medbreak\item   [$\circ$] \cite{Gromov} If  $\varGamma$ is a finitely generated group, then
$\GK \ku \varGamma < \infty$  if and only if  $\varGamma$  is nilpotent-by-finite. 
\end{itemize}

\subsection{Noetherian Hopf algebras}\label{subsec:Noetherian}

A Hopf algebra is  \emph{Noetherian}, respectively \emph{Artinian}, if the underlying algebra is so, 
that is, it satisfies the ascending, respectively descending,  chain condition on left ideals. 
Artinian Hopf algebras are finite-dimensional \cite{Liu-Zhang}, so we focus on the Noetherian ones.
The following result is classical.

\begin{lemma}\label{lemma:filtration-noeth}  
Let $A$ be a filtered algebra. If $\gr A$ is Noetherian, then $ A$ is Noetherian.
\end{lemma}

\begin{itemize}[leftmargin=*]
 \medbreak\item   [$\circ$] If a Lie algebra $\g$ is finite-dimensional, then  $U(\g)$ is Noetherian; 
 the converse is the well-known open Question \ref{question:enveloping-Noetherian}. \phantomsection\label{text-question:enveloping-Noetherian} A recent relevant result is the following:
\end{itemize}

\begin{theorem} \cite{Sierra-Walton}
The  enveloping algebra of the Witt algebra is not Noetherian. \qed
\end{theorem}

\begin{itemize}[leftmargin=*]
\medbreak\item   [$\circ$] The group algebra of a polycyclic-by-finite group is Noetherian \cite{Hall}; the converse is
the well-known open Question \ref{question:group-algebra-Noetherian}. \phantomsection\label{text-question:group-algebra-Noetherian}

\medbreak
\item   A recent relevant result: If $\ku \varGamma$ is Noetherian,
then $\varGamma$ is amenable and every subgroup of $\varGamma$ is finitely generated; see
\cite[Corollary B]{Kropholler-Lorensen}, based on  \cite{Bartholdi}.
\end{itemize}

\begin{remark} 
A group is Noetherian if it satisfies the maximal condition on subgroups. 
If the group algebra of a given group is Noetherian, then so is the group but the converse is not true.
Indeed there are  Noetherian groups that are not polycyclic-by-finite \cite{Olshanski} but their group algebras
are not Noetherian \cite{Ivanov}.  
\end{remark}

In general an object in a category
is Noetherian if it satisfies the maximal condition on subobjects. 
Say as in \cite{Drinfeld} that the category \QG{\ } of quantum groups is the
opposite to the category \Hopf{\ } of Hopf algebras with bijective antipode. Then we  may consider  
\emph{Noetherian quantum groups} (Noetherian objects in \QG) 
and  Noetherian objects in \Hopf.
To my knowledge these have not been studied.

\medbreak
Now a linear Noetherian group  is polycyclic-by-finite, as  follows from a
remarkable result of Tits:

\begin{theorem}\label{thm:linear-noeth} \cite{Tits}
A finitely generated linear group is either solvable-by-finite or contains a non-cyclic free subgroup.\qed
\end{theorem}

Complementing this result we mention: 
\begin{theorem}\label{thm:polycyclic-linear} (Auslander, Swan)
Any polycyclic-by-finite group is linear.  \qed
\end{theorem}

\begin{theorem} 
 \cite{Molnar}    A commutative Hopf algebra $H$  is  Noetherian iff it is affine. \qed
\end{theorem}

If this is the case, then $H$ has  finite $\GK$, which equals the Krull dimension, but the converse is false. 
Question \ref{question:Noetherian-implies-affine} asks whether the commutativity in the preceding result is needed.
\phantomsection\label{text-question:Noetherian-implies-affine}

\begin{remark} \label{text-question:affine+finiteGK-Noetherian}
If $G$ is not finitely generated, but any finitely generated subgroup is finite (e.g. $\Gb_{\infty}$), then
$\GK \ku G =0$ but $\ku G$ is not Noetherian.  However it is plausible that 
``$\GK H < \infty \implies$ Noetherian'' holds when $H$ is affine, see
 Question \ref{question:affine+finiteGK-Noetherian}.
\end{remark}

 Lemma \ref{lemma:filtration-noeth} also gives:
\begin{itemize}[leftmargin=*]
\medbreak\item   [$\circ$]  If $\g$ is a semisimple finite-dimensional Lie algebra, 
then $U_q(\g)$ is Noetherian.
\end{itemize}

\begin{remark} \label{rem:text-noetherian}
Two   facts that make  difficult the setting up of a classification strategy:  

\begin{itemize}[leftmargin=*]
\medbreak\item   [$\circ$]   The converse in  Lemma \ref{lemma:filtration-noeth}
does not hold.

\medbreak\item   [$\circ$]   A subalgebra of a Noetherian algebra is not necessarily Noetherian. 
One wonders what happens for Hopf algebras, see Question \ref{question:subalg-Noetherian}.
 \label{text-question:subalg-Noetherian}
\end{itemize}
\end{remark}

\subsection{Co-Frobenius and cosemisimple Hopf algebras}\label{subsec:Cofrob}

Let $H$ be a Hopf algebra. We start with two definitions.

\begin{itemize}[leftmargin=*]
\medbreak\item     \cite{Sweedler-integral} A (left)  \emph{integral} on $H$ is a linear functional $ \smallint \in H^*$ 
invariant under   the left regular representation; that is,
$\langle \smallint, x\rangle =  x\_{1} \langle \smallint, x\_{2}\rangle$,  for all $x  \in H$.
The map $\smallint$ is reminiscent of the Haar measure of a Lie group, hence the name. 
The set $\Ic_{\ell}(H)$ of all left integrals on $H$ is a subspace of $H^*$
and $\dim \Ic_{\ell}(H) \leq 1$; we say that $H$ is \emph{co-Frobenius} iff $\dim \Ic_{\ell}(H) = 1$, i.e.,
there exists a non-zero left integral on $H$. 
Analogously $ \smallint \in H^*$ is a right integral if
$\langle \smallint, x\rangle =  \langle \smallint, x\_{1}\rangle x\_{2} $,  for all $x  \in H$.
The subspace $\Ic_r(H)$ of all right integrals satisfies $\dim \Ic_r(H) = \dim \Ic_\ell(H) $.

\medbreak\item    $H$ is (left) cosemisimple   iff   the category $\lcomod{H}$ of left $H$-comodules is semisimple.
\end{itemize}

Let us discuss some properties and examples.

\begin{itemize}[leftmargin=*]
\medbreak\item   [$\circ$]  \cite{Sweedler-integral}
$H$ is left cosemisimple iff $H \simeq \bigoplus\limits_{V \in \lcomod{H}\text{ simple}} C_V$, 
iff $H$ is right cosemisimple, 
iff $ H$ is co-Frobenius and for some (hence any) $\smallint \in \Ic_{\ell}(H)$ one has $\smallint (1) \neq 0$.

\medbreak\item   [$\circ$] \cite{Larson-Sweedler} If $\dim H < \infty$, then $H$ is co-Frobenius (but not always cosemisimple).

\medbreak\item   [$\circ$] \cite{Larson-Sweedler}
The following data are equivalent:
(a) a right integral $\smallint: H \to \ku$, 
(b) a bilinear form $(\, \vert \, ): H \times H \to \ku$ satisfying
\begin{align} \label{eq:killing1}
( xy\vert z ) &= ( x\vert yz ), & x,y,z &\in H,
\\\label{eq:killing2}
( f \rightharpoonup x\vert y ) &= (  x\vert \Ss(f) \rightharpoonup y), & x,y &\in H,
f\in H^*.
\end{align}
\end{itemize}
Explicitly, $\langle \smallint, x\rangle = (x \vert 1)$, $(x \vert y) = \langle \smallint, xy\rangle$, for any
$x,y\in H$. When $H$ is cosemisimple and $\langle \smallint, 1\rangle = 1$, the corresponding
$(\, \vert \, )$, called the \emph{Killing form} of $H$, is non-degenerate.

\begin{itemize}[leftmargin=*] 
\medbreak\item   [$\circ$]  The enveloping algebra of a Lie algebra is never co-Frobenius.

\medbreak\item   [$\circ$]  Every group algebra  is cosemisimple.

\medbreak\item   [$\circ$]  \cite{Sullivan} Let $\Gb$ be an algebraic group. Then  
$\Oc(\Gb)$ is cosemisimple   iff  $\Oc(\Gb)$  is co-Frobenius   iff  $\Gb$ is reductive.

\medbreak\item   [$\circ$]  If $\g$ is a simple Lie algebra, then
$U_q(\g)$ is never co-Frobenius.

\medbreak\item   [$\circ$]  Let $\Gb$ be  a simple algebraic group. Then 
$\Oc_q(\Gb)$ is co-Frobenius, but $\Oc_q(\Gb)$ is cosemisimple    iff  $q$ is not a root of 1.
\end{itemize}

Below we  summarize a series of characterizations of co-Frobenius Hopf algebras
that enhance its interest, due to several people:
S.~D\u{a}sc\u{a}lescu, C.~N\u{a}st\u{a}sescu, S.~Donkin, B.~I-Peng Lin, D. E. Radford and the 
authors of \cite{A-Cuadra,A-Cuadra-Etingof}--see loc. cit. for references.
Given $M \in \lcomod{H}$, its injective hull is denoted by $E(M)$. 

\begin{theorem}\label{th:oldchar} 
The following statements are equivalent:
\begin{multicols}{2}
\begin{enumerate} [label=\rm{(\roman*)}] 
\medbreak\item    $H$ is co-Frobenius.
\medbreak\item   \label{item:findim-injhull} $\dim E(S) < \infty$   $\forall\,S \in \lcomod{H}$ simple.
\medbreak\item    $E(\ku)$ is finite-dimensional.
\medbreak\item    There exists $ Q \in\lcomod{H}$ injective, $ \dim Q < \infty$, $Q \neq 0$.
\medbreak\item    Every $0\neq M \in \lcomod{H}$ has a  finite-dimensional quotient $\neq 0$.
\medbreak\item    $\lcomod{H}$ has a nonzero projective.
\medbreak\item    Every $M \in \lcomod{H}$ has a projective cover.
\medbreak\item    Every injective in $\lcomod{H}$ is projective.
\medbreak\item    The right versions of  these characterizations.
\medbreak\item    The coradical filtration of $H$ is finite. \qed
\end{enumerate}
\end{multicols}
\end{theorem}

Question \ref{question:cofrob-coss} is natural in the spirit of this article. 
\phantomsection\label{text-question:cofrob-coss} We point out:

\begin{itemize}[leftmargin=*] 
\medbreak\item   [$\circ$]  The crucial step towards the classification of pointed co-Frobenius
Hopf algebras with abelian group is given in \cite[Theorem 1.4.2]{A-Angiono-Heckenberger-triang}, using  results from \cite{A-Dasca,A-Cuadra-Etingof}. Namely, the corresponding
 diagram should be a finite-dimensional  Nichols algebra of diagonal type. See Section \ref{section:pointed} for unexplained terminology.

\medbreak\item   [$\circ$]  It seems that  co-Frobenius Hopf algebras could be built-up from cosemisimple and  finite-dimensional ones. See \cite{A-Cuadra-Etingof,ANT} for 
examples disproving a more precise conjecture from \cite{A-Dasca} but still supporting this belief.

\medbreak\item   [$\circ$]  The classification of the finite-dimensional cosemisimple (hence semisimple)
Hopf algebras is known only in  low dimensions. See \cite{Natale} for a recent account.
\end{itemize}

We continue the discussion on cosemisimple Hopf algebras in Subsection \ref{subsec:compact}.

\subsection{The antipode}\label{subsec:Antipode}
By convenience we assume that the antipode of any Hopf algebra  is bijective.
There are examples where this does not hold (Takeuchi, Nichols, see \cite{Schauenburg-non-fflat}) 
but, intuitively, they \emph{grow too fast}.
Skryabin proved that the antipode of a Noetherian Hopf algebra is
injective and conjectured that it is bijective \cite{Skryabin-antipode}, 
cf. Conjecture \ref{conjecture:antipode}. 
The next result follows from \cite[Theorem A]{Skryabin-antipode} as kindly explained to me by K. Goodearl.
\begin{theorem} 
The antipode of a Hopf algebra  domain with finite  Gelfand-Kirillov dimension is bijective. \qed
\end{theorem}

Here is   another indication of the validity of the Conjecture, see footnote \ref{foot:res-fd}.

\begin{theorem} \cite{Skryabin-antipode-Noetherian}
The antipode of an either right or left Noetherian residually finite-dimensional Hopf algebra  is bijective.  \qed
\end{theorem}

The antipode of a finite-dimensional Hopf algebra is bijective \cite{Larson-Sweedler}.
This was extended to co-Frobenius Hopf algebras by Radford \cite{Radford-cofrob}.
There are results with homological flavor: the antipode of a  Hopf algebra that has Van den Bergh duality is bijective
\cite{Le Meur}; the same under some finite homological conditions \cite{Lu-et-al}.

\subsection{Regular Hopf algebras}\label{subsec:Regular}

A Noetherian algebra is regular if it has finite global dimension; see \cite[Chapter 7]{Mcconell-Robson}.  
The following Noetherian Hopf algebras are regular:
\begin{itemize}[leftmargin=*]
\item [$\circ$]   $U(\g)$, where $\g$ is a finite-dimensional Lie algebra.

\medbreak\item[$\circ$]    $\ku \varGamma$, where $\varGamma$ is polycyclic-by-finite \cite[7.5.6]{Mcconell-Robson}.

\medbreak\item[$\circ$]    $\Oc(\Gb)$, where $\Gb$  is an algebraic group. 
\end{itemize}

It is well-known that a finite-dimensional Hopf algebra, being Frobenius, either is  
semisimple, or has infinite global dimension. The following useful criterion follows from this,
 together with a theorem of Schneider.

\begin{lemma}\cite{ANT}
A Hopf algebra having a finite-dimensional  non-semisimple normal Hopf subalgebra has infinite global dimension.
\end{lemma}
Thus Lusztig's  $\mathcal U_\epsilon(\g)$ is not regular.
For a more in-depth discussion of the homological properties of Hopf algebras see the surveys \cite{Brown-survey,Brown-Gilmartin-survey,Brown-Zhang-survey,Goodearl-survey}.

\subsection{Semiprimitive Hopf algebras and the Jacobson radical}\label{subsec:Semiprimitive}
Let $A$ be an algebra. It is well-known that  its Jacobson radical  is 
\begin{align*}
\Jac A \coloneqq \bigcap_{V \in \Irr A} \ker \rho_V = \bigcap_{\Mg \text{ maximal left ideal}} \Mg
= \bigcap_{\Ng \text{ maximal right ideal}} \Ng =
\end{align*}
the  largest ideal $I$ such that $1- x$ is a unit for each $x\in I$.
Recall that  $A$ is  
\begin{enumerate}[leftmargin=*,label=\textrm{(\roman*)}]
\medbreak\item   \label{item:primitive} \emph{primitive} if there exists $V \in \Irr A$ with $\ker \rho_V = 0$,
\medbreak\item   \label{item:semiprimitive} \emph{semiprimitive} if  $\Jac A = 0$.
\end{enumerate}
Clearly \ref{item:primitive} $\implies$ \ref{item:semiprimitive}. See \cite{Passman-survey} for the history of the next result.
\begin{itemize}[leftmargin=*]
\medbreak\item[$\circ$]     Every group algebra  is semiprimitive provided that $\ku \neq  \overline{\Q}$ (when $\ku =  \overline{\Q}$
it is open; when $\car \ku > 0 $
the situation is much more delicate, see \cite{Passman-survey}).

\medbreak\item[$\circ$]  There exist group algebras  that are  primitive \cite[Chapter 9]{Passman-book}.

\medbreak\item[$\circ$]  For any finite-dimensional Lie algebra $\g$, $U(\g)$ is semiprimitive \cite[3.1.16]{dixmier}.

\medbreak\item[$\circ$]    The finite-dimensional Lie algebras $\g$ with $U(\g)$  primitive are
described in \cite{Ooms}.
\end{itemize}

\subsection{Hopf algebras that are domains, prime, or  semiprime}\label{subsec:Domain}

Recall that an algebra $A$ is
\begin{enumerate}[leftmargin=*,label=\textrm{(\roman*)}] \setcounter{enumi}{2}
\medbreak\item   \label{item:domain}  a \emph{domain} if $xy = 0$ implies $x = 0$ or $y = 0$,
\medbreak\item   \label{item:prime} \emph{prime} if  $xAy = 0$ implies $x = 0$ or $y = 0$,
\medbreak\item   \label{item:semiprime} \emph{semiprime} if $xAx = 0$ implies $x = 0$,
\end{enumerate}
for every $x, y \in A$. Clearly one has \hspace{5pt}
$\begin{aligned}
\text{\ref{item:domain}} \implies   & \text{ \ref{item:prime}} \implies  \text{\ref{item:semiprime}}
\\
 & \,\Uparrow  \qquad\quad \, \, \Uparrow
\\ 
&\text{ \ref{item:primitive}}  \implies \, \text{\ref{item:semiprimitive}}
\end{aligned}$,
 see \cite{Mcconell-Robson} for details.
  Again the next statement is classical.

\begin{lemma}\label{lemma:filtration-domain}  
Let $A$ be a filtered algebra. If $\gr A$ is domain, then so is $A$.
\end{lemma}
\begin{itemize}[leftmargin=*]
\medbreak\item[$\circ$]      The enveloping algebra of any Lie algebra is a domain.

\medbreak\item[$\circ$]    The enveloping algebra of a finite dimensional Lie superalgebra $\lgo$ 
with $\lgo_1 \neq 0$ is not a domain; it is seldom prime and it could fail to be semiprime.
The simple $\lgo$ with $U(\lgo)$ prime by the criterium in \cite{Bell} are listed
in  \cite{Wilson}; the case of  $\lgo =W(n)$ for odd $n > 5$ is still open.
\end{itemize}

Let $\varGamma$ be a group. 
If $\ku \varGamma$ is a domain, then $\varGamma$ is torsion-free: if $1 \neq x\in \varGamma$
has finite order, then $1-x$  is a zero divisor.
The converse, Conjecture \ref{conjecture:0-divisor}, 
is the famous zero-divisor conjecture of Kaplansky, which is still open. \phantomsection\label{text:conjecture:0-divisor}
We record a case where it is known to be true.

\begin{theorem} \label{thm:kropholler} \cite{Kropholler-et-al}
If $\varGamma$ is a torsion-free solvable group, then $\ku\varGamma$  is a domain. \qed
\end{theorem}

 \begin{itemize}
	\item  [$\circ$]  Every group algebra  is semiprime \cite[Theorem 4.2.13]{Passman-book}.
\end{itemize}

Prime group algebras are characterized by the following result:
\begin{theorem}\label{thm:prime-group-algebras} \cite[Theorem 8]{Connell} The following are equivalent:
\begin{enumerate}[leftmargin=*,label=\textrm{(\roman*)}]
\medbreak\item   \label{item:prime-kuG} The algebra $\ku\varGamma$ is prime.

\medbreak\item    $FC(\varGamma)$ is a torsion-free abelian subgroup.

\medbreak\item    If $N \trianglelefteq \varGamma$ is a finite normal subgroup, then $N$ is trivial.   \qed
\end{enumerate}
\end{theorem}

\begin{corollary}\label{cor:prime-nilpotent}
If $\varGamma$ is  finitely generated  nilpotent, then
  the following  are equivalent:
\begin{enumerate}[leftmargin=*,label=\textrm{(\roman*)}]
\medbreak\item   \label{item:torsion-free-G-nilp}  $\varGamma$ is torsion-free.
\medbreak\item   \label{item:domain-kuG-nilp}  The algebra $ \ku \varGamma$ is a domain.
\medbreak\item   \label{item:prime-kuG-nilp}  The algebra $ \ku \varGamma$ is prime.
\end{enumerate}
\end{corollary}

\pf   \ref{item:torsion-free-G-nilp} $\implies$ \ref{item:domain-kuG-nilp}: by Theorem \ref{thm:kropholler}.
\ref{item:domain-kuG-nilp} $\implies$ \ref{item:prime-kuG-nilp}: is clear.
\ref{item:prime-kuG-nilp} $\implies$ \ref{item:torsion-free-G-nilp}: it is well-known that the torsion elements of $G$
form a subgroup $T$ of $G$ which is finite and normal. 
Then Theorem \ref{thm:prime-group-algebras} implies that $T$ is trivial. 
\epf

The natural Question \ref{question:classification:domain-prime-semiprime},
 seems to be too ambitious unless approached
under various restrictions.
\phantomsection\label{text-question:classification:domain-prime-semiprime}
For instance, a domain is not always semiprimitive, 
but  no example among Hopf algebras seems to be known,
 raising Question \ref{question:domain-semiprimitive}. 
\phantomsection\label{text-question:domain-semiprimitive}

\section{Pointed Hopf algebras}\label{section:pointed}

\subsection{The lifting method}\label{subsec:Method}

Let $H$ be a Hopf algebra such that its coradical $H_0$ is a Hopf subalgebra of $H$. Then 
the coradical filtration $(H_n)_{n\in \N_0}$ is a Hopf algebra filtration \cite[5.2.8]{Montgomery-book}
and the corresponding graded coalgebra $\gr H$ 
is in fact a Hopf algebra that splits as the  bosonization
\begin{align*}
\gr H \simeq \Rc \# H_0
\end{align*}
where $\Rc = \oplus_{n\ge 0} \Rc^n$ is a graded connected Hopf algebra in $\yd{H_0}$.
In fact, $\gr H$ is coradically graded (the coradical filtration coincides with the canonical filtration arising from the grading) hence so is $\Rc$. Let $V = \Rc^1$; $\Rc$ is called the \emph{diagram}, and $V$ the
\emph{infinitesimal braiding}, of $H$.
Then the subalgebra of $\Rc$ generated by $V$ is isomorphic to the Nichols algebra $\toba(V)$, 
see below. 
These observations motivate the following method, proposed in \cite{A-Schneider-jalg-p3,A-Schneider-adv,A-Schneider-cambr}, 
for the classification of Hopf algebras
whose coradical is a Hopf subalgebra, and have either finite dimension,
or finite $\GK$, or any suitable property; but observe that as of today,
the method applies poorly to Noetherian Hopf algebras, 
see Remark \ref{rem:text-noetherian}.
It consists of several steps:

\begin{enumerate} \setcounter{enumi}{-1}
\medbreak\item  \label{item:levante0} Choose an  adequate cosemisimple Hopf algebra $K$ and describe $\yd{K}$.
\end{enumerate}

\begin{enumerate}[leftmargin=*,label=\textrm{(\alph*)}]
\medbreak\item \label{item:levante1} Classify those $V\in \yd{K}$ such that $\dim \toba(V) < \infty$ or $\GK \toba(V) < \infty$.

\medbreak\item  \label{item:levante2} Given $V$ as in \ref{item:levante1}, determine all coradically graded connected Hopf algebras
$\Rc \in \yd{K}$ such that $\Rc^1 \simeq V$ and $\dim \Rc < \infty$ or $\GK \Rc < \infty$.

\medbreak\item \label{item:levante3} Given $\Rc$ as in \ref{item:levante2}, find all Hopf algebras $H$ such that
$\gr H \simeq \Bc \# K$.
\end{enumerate}

See \cite{A-Schneider-crelle,A-Schneider-annals} for some achievements within this method.
We now discuss some relevant features  referring for more details to later parts of this text.
Note that we assume most of the time that $K$ and $\Rc$ are affine, which implies that $\gr H$ and, a fortiori, $H$
are also  affine; but there braided vector spaces $V$ arising from abelian groups  
with $\dim V = \infty$ and $\GK \toba(V) < \infty$, see \cite{A-Angiono-Heckenberger-infinite}. 
The analysis in loc. cit. leads to Question \ref{question:nichols-locally-finite}. \phantomsection\label{text-question:nichols-locally-finite}
It seems to be open whether $H$ affine implies that $\gr H$ is affine, Question \ref{question:H-affine-grH-affine}. \phantomsection\label{text-question:H-affine-grH-affine}

\medbreak
\noindent \emph{Step \ref{item:levante0}}: Here by adequate we mean related to the problem at hand;
e.g. in the finite-dimensional context one has to assume that $\dim K < \infty$. 
 Then $K$ is also semisimple (and vice versa, by \cite{Larson-Radford}) 
 and a fortiori $\yd{K}$ is a semisimple category.  
However, when turning to finite $\GK$,  if $K$ is infinite-dimensional, then it is not semisimple and $\yd{K}$ 
is mostly wild.
Still, the condition $\GK \toba(V) < \infty$ is very restrictive and at least in some substantial cases 
a rather complete answer of \ref{item:levante1} is envi\-sageable.
Now cosemisimple Hopf algebras with finite $\GK$ are not yet understood\footnote{beyond the classical and quantum cases, see Section \ref{section:more}}
 so we focus on the simplest case:

\begin{definition}
A Hopf algebra $H$ is \emph{pointed} if its coradical is a group algebra, i.e. $H_0 \simeq  \ku \varGamma$,
where $\varGamma \simeq G(H)$.
\end{definition} 

Thus, to classify finite-dimensional pointed Hopf algebras one has to go through the steps  \ref{item:levante1},
\ref{item:levante2} and \ref{item:levante3} under the assumption $K \simeq \ku \varGamma$, where $\varGamma$
is a finite group, while the analogous problem with finite $\GK$ requires  $\varGamma$ to be nilpotent-by-finite.
If we look at pointed Hopf algebras that are domains and have finite $\GK$, then we need $\varGamma$ to be nilpotent-by-finite and torsion-free, see Corollary \ref{cor:prime-nilpotent}.

\begin{remark} \cite{A-Galindo-Muller}  As said, the classification of the semisimple Hopf algebras is presently out of
reach but for some of them the problem is approachable as we argue next.
Two semisimple Hopf algebras $L$ and $K$ are Morita-equivalent  iff there is an equivalence $\Gc:\yd{L} \to \yd{LK}$ of braided tensor categories (this is not the same as being Morita-equivalent as algebras).
If this happens, then $\Gc$
preserves Nichols algebras, i.e., $\Gc(\toba(V)) \simeq \toba(\Gc(V ))$, for $V \in \yd{L}$.  
Thus is we know some, or all, finite GK-dimensional Nichols algebras in $\yd{L}$ (e.g. when $L\simeq \ku \varGamma$
for specific $\varGamma$), then this knowledge is transfered to $\yd{K}$. Also, 
when there are no finite GK-dimensional Nichols algebras in $\yd{L}$, the same holds for $K$.
\end{remark} 

\medbreak
\noindent \emph{Step \ref{item:levante1}}:
Given $V \in \yd{K}$, where $K$ is any Hopf algebra, there exists a unique graded connected Hopf algebra  
$\toba = \oplus_{n\ge 0} \toba^n$ in $\yd{K}$
which is generated by $\toba^1 \simeq V$ and is coradically graded. It is called the Nichols algebra of $V$ and denoted by
$\toba(V)$.
 See   \cite{A-Angiono-Diag-survey,A-Schneider-cambr,A-Leyva,Heckenberger-Schneider-book} for detailed expositions.  
Here we just stress:

\begin{itemize}[leftmargin=*] 
\medbreak\item  The Nichols algebra $\toba(V)$ depends, as algebra and coalgebra, only on the underlying braided vector space  \eqref{eq:braid-yd}; thus it is more efficient to look at Nichols algebras of braided vector spaces. 
Besides, no all-embracing approach to Step \ref{item:levante1} is known (yet) so 
we proceed by looking at suitable classes of braided vector spaces.

\medbreak\item  The  braided vector spaces of diagonal type are the best understood. Recall that $(V,c)$ is of diagonal type
if $V$ has a basis $(v_i)_{i\in I}$ such that $c(v_i \otimes v_j) = q_{ij}v_j \otimes v_i$ for all $i, j\in I$, where the
$q_{ij}$'s belong to $\kut$.
We summarize and refer to \cite{A-Angiono-Diag-survey}:
\end{itemize} 

\begin{itemize}[leftmargin=*]
\medbreak\item[$\circ$]  The braided vector spaces $(V,c)$ of diagonal type with 
$\dim \toba(V) < \infty$ were classified in \cite{Heckenberger-classif RS} by means of the Weyl groupoid. 
Actually the paper gives the classification of the $(V,c)$ of diagonal type with `arithmetic root system'
see e.g. \cite{A-Angiono-Diag-survey,Heckenberger-Schneider-book} for details.
The classification was organized in \cite{A-Angiono-Diag-survey} in four types:

\begin{multicols}{2}
\begin{enumerate}[leftmargin=*,label=\textrm{(\alph*)}]
\medbreak\item  Cartan type  \cite{A-Schneider-adv};
\medbreak\item  super type \cite{Heckenberger-Yamane};
\medbreak\item  (super) modular type;
\medbreak\item  UFO's (the rest).
\end{enumerate}
\end{multicols}

\medbreak\item[$\circ$]  
The defining relations of the Nichols algebras with arithmetic root system were given in \cite{Angiono-crelle}.

\end{itemize} 

\medbreak
When $\varGamma$ is a finite abelian group, any $V\in \yd{\ku \varGamma}$  is of diagonal type. 
In the context of finite $\GK$, when $\varGamma$ is infinite, 
this no longer true; see \S \ref{subsec:Nichols} for a discussion.

\begin{itemize}[leftmargin=*]
\medbreak\item[$\circ$]   The  next class of braided vector spaces  that received considerable attention is that coming from 
non-abelian finite groups. See \cite{A-Leyva,A-Carnovale-Garcia-VII} for details and references.
Succintly, the decomposable braided vector spaces over finite groups with finite-dimensional Nichols algebras are classified \cite{Heckenberger-Vendramin-rank>2}. 
For the indecomposable ones there is a small set of examples with $\dim\toba(V) < \infty$, a few criteria insuring infinite dimension and various families that are completely open. 
We discuss finite GK-dimensional pointed Hopf algebras over non-abelian groups
in \S \ref{subsec:Nichols-nonabelian}.
\end{itemize}

\medbreak
\noindent \emph{Step \ref{item:levante2}}: 
Let $V$ be a braided vector space.
For brevity, we say that a coradically  graded connected Hopf algebra 
$\Rc = \oplus_{n\ge 0} \Rc^n$  such that $\Rc^1 \simeq V$ is a  \emph{post-Nichols algebra} of $V$. \phantomsection \label{text-defn:post-nichols}
Let $K$ be a finite-dimensional cosemisimple Hopf algebra, let  
$V\in \yd{K}$ such that $\dim \toba(V) < \infty$
and let $\Rc \in \yd{K}$ be a post-Nichols algebra of $V$
such that  $\dim \Rc < \infty$.
It was conjectured that $\Rc \simeq \toba(V)$, cf. \cite[1.4]{A-Schneider-adv} for $K\simeq \ku \varGamma$
and \cite[2.7]{A-Schneider-cambr} for general $K$. The Conjecture was proved for $V$ of Cartan diagonal type
in \cite{A-Schneider-annals} and for general diagonal type  in \cite{Angiono-crelle}. It was also checked for 
the known examples arising from finite non-abelian groups.

\medbreak
In the context of finite $\GK$, the conjecture is no longer true. Our approach
to this question is summarized in \S \ref{subsec:pre-post-nichols}. 
As an illustration, we discuss post-Nichols algebras of vector spaces in \S \ref{subsec:post-Nichols-vs}.

\medbreak
\noindent \emph{Step \ref{item:levante3}}: This is a deformation question: given $\gr H \simeq \Rc \# K$ 
compute all possible $H$ (called liftings).
Remarkably, in the finite-dimensional context, 
this is rigid enough to be approachable by combinatorial means, giving rise to explicit
answers in many cases \cite{A-Schneider-annals,Angiono-Garcia}. 
These techniques could be adapted to finite $\GK$, but there are crucial new features:

\begin{itemize}[leftmargin=*]
\medbreak\item  In the finite-dimensional setting, being $\ku \varGamma$  semisimple, we have a decomposition 
\begin{align}\label{eq:decomposition-first-corad}
H_1 \simeq H_1/H_0 \oplus H_0,
\end{align}
\end{itemize} 
of $\varGamma$-modules; this is no longer true for liftings  
over infinite groups, see  Example \ref{exa:corrigendum}.

\begin{itemize}[leftmargin=*]
\medbreak\item  In particular, the infinitesimal braiding $V$ of $H$ may contain a subspace $V'$ with basis $(v_i)$ such 
that $c(v_i \otimes v_j) = v_j \otimes v_i$ for all $i,j$; thus Lie algebras enter into the picture as deformations 
of $\toba(V') \simeq S(V')$. See Example \ref{exa:double-jordan}.

\medbreak\item    In the finite-dimensional setting, $\Rc$ is a Nichols algebra (conjecturally, see above) but in general
liftings of post-Nichols algebras have to be determined too. 
\end{itemize}

\begin{example}\label{exa:corrigendum} \cite{A-Angiono-Heckenberger-Jordan-corrig}
See \eqref{eq:jordan-plane} for the Jordan plane.
Let $\mathfrak U$ be the Hopf algebra generated by $a_1, a_2, g, g^{-1} $ 
with defining relations and comultiplication given by
\begin{align*}
ga_1 &= a_1g + (g-g^2), &  ga_2 &= a_2g + a_1g, & g^{\pm1}  g^{\mp1} &= 1, 
\\  a_1a_2 - a_2a_1 &= \dfrac{a_1^2}{2} - a_2 -\dfrac{1}{2}a_1, 
& \Delta (a_i) &= a_i \otimes g + 1 \otimes a_i, \ i= 1,2,& \Delta (g) &= g \otimes g.
\end{align*}
Then $\mathfrak U$ is pointed, 
$\gr \mathfrak U \simeq J \# \ku \Z$, where $J$ is the Jordan plane, but \eqref{eq:decomposition-first-corad} does not hold.
\end{example}

\begin{example}\label{exa:double-jordan}
The double of the Jordan plane is the Hopf algebra $\D$  introduced in \cite{A-Penha},
presented by generators $u$, $v$, $\xi$, $g^{\pm 1}$, $x$, $y$ and defining relations
\begin{align*}
&\begin{aligned}
g^{\pm 1} g^{\mp 1}&= 1, &\xi  g &= g \xi, &yx &= xy - \frac{1}{2} x^2, & vu &= uv - \frac{1}{2}u^2,
\\
gx &= xg, & gy &= yg +xg, & \xi  y &= y \xi  -2 y, &  \xi  x&=x \xi  - 2 x, 
\\
u g &= gu, & v g &=gv  + gu, & v\xi  &= \xi v  -2 v, &  u\xi  &= \xi u -2 u,
\end{aligned}
\\
&\begin{aligned}
u x &=x u , &   v x &= x v  + (1-g) + xu, &
u y &=yu  +(1 - g), &  v y &=y v + \frac{1}{2}g \xi  + y u. 
\end{aligned}
\end{align*}
The  Hopf algebra structure is given  by 
$g\in G(\D)$, $u,\xi \in\Prim(\D)$, $x,y\in\Prim_{g,1}(\D)$ and 
\begin{align*}
\Delta( v) &=  v\otimes 1 + 1\otimes v -\frac{1}{2}\xi\otimes u. 
\end{align*}
It is an affine Noetherian domain of $\GK = 6$; it is pointed,
the infinitesimal braiding $V$ is generated by $x,y,u,\xi$ while $v$ belongs to the
diagram $\Rc$ but not to $\toba(V)$. Note that $\ku\langle u,\xi\rangle \simeq U(\mathfrak b_2)$ where $\mathfrak b_2$ is the  Borel subalgebra of $\mathfrak{sl}(2)$.
\end{example}

\subsection{Nichols algebras over abelian groups}\label{subsec:Nichols}
Let $\varGamma$ be an abelian group and $V \in \yd{\ku \varGamma}$, $\dim V < \infty$.
Now the Yetter-Drinfeld modules over $\ku\varGamma$ are the $\varGamma$-graded $\varGamma$-modules;
set $V= \oplus_{g\in \varGamma} V_g$ for the $\varGamma$-grading. If $V^{(\chi)} $ is as in \eqref{eq:decomp-commutative} for $\chi \in \widehat{\varGamma} \simeq \Alg(\ku \varGamma, \ku)$,
and $V_g^{(\chi)} = V_g \cap V^{(\chi)}$, then \eqref{eq:decomp-commutative-dix} says that
\begin{align*}
V &= \bigoplus_{\substack{g\in \varGamma,  \chi \in \widehat{\varGamma}}} V_g^{(\chi)}.
\end{align*}
 
 Thus, given a decomposition 
$V = \oplus_{i\in \I_{\theta}} V_{i}$ in $\yd{\ku\varGamma}$, with the $V_{i}$'s indecomposable, there exist 
$g_i \in \varGamma$ and $\chi_i \in \widehat{\varGamma}$ such that $V_{i} \subseteq V_{g_i}^{(\chi_i)}$ for $i\in \I_{\theta}$.
Two invariants  are:
\begin{itemize}[leftmargin=*]
\medbreak\item  the matrix $\bqg = (q_{ij})_{i,j \in \I_{\theta}}$,  $q_{ij} = \chi_j(g_i)$;

\medbreak\item  the family $\ag = (\ag_{ij}) _{i,j\in \I_{\theta}}$, $\ag_{ij}\in \End V_j$, determined by
${g_i}_{\vert V_j} = q_{ij}\big(\id_{\vert V_j} + \ag_{ij} \big)$. Clearly,  
$\ag_{ij}$ is nilpotent and
$\ag_{ik}\ag_{jk} = \ag_{jk}\ag_{ik}$ for all $i,j,k$; 
a suitable normalization of $\ag$ is called the \emph{ghost}.
Then the braiding of $V$ is given by
\begin{align}\label{eq:braiding-abelian}
c(x\otimes y) &= q_{ij}\left(y + \ag_{ij}(y) \right) \otimes x, & x\in V_i, y\in V_j.
\end{align}
\end{itemize}

We distinguish three kinds of indecomposables:

\begin{itemize}[leftmargin=*]
\medbreak\item    $V_i$ is  a \emph{point} if it is simple, i.e., has dimension 1 and braiding given by $q \in \kut$;

\medbreak\item  $V_i$ is  a \emph{block} if it is indecomposable as braided vector space but   
 $\dim V >1$;

\medbreak\item   $V_i$ is  a \emph{pale block} if it is indecomposable in $\yd{\ku \varGamma}$ but not
as braided vector space.
\end{itemize}
When $V$ is semisimple, it is  of diagonal type and the only invariant is the braiding matrix $\bqg$.
In this case we set $\toba_{\bqg} \coloneqq \toba(V)$ and represent the matrix $\bqg$ (up to twist-equivalence)
by its Dynkin diagram, that has $\theta$ vertices with the $i$-th  vertex labelled by $q_{ii}$;
 there is one edge between $i$ and $j$, labeled by $\widetilde{q}_{ij} \coloneqq q_{ij}q_{ji}$,
when $\widetilde{q}_{ij} \neq 1$, otherwise there is none. 
That is, $\xymatrix{\cdots \underset{ i }{\overset{q_{ii}}{\circ}} \ar  @{-}[r]^{ \widetilde{q}_{ij}}  &  \underset{ j }{\overset{q_{jj}}{\circ} } \cdots}$.  A connected component of $\bqg$ is
the submatrix spanned by a connected component of its Dynkin diagram.

\medbreak
We summarize the state of the classification of those $V$ with $\GK \toba(V) < \infty$.

\begin{theorem}\label{thm:abelian-case}
\emph{(i)} \cite{Angiono-Garcia2} If $\GK \toba(V) < \infty$ and $V$ is semisimple, i.e.,  of diagonal type, then $V$ has arithmetic root system, hence its connected components  appear in the list of \cite{Heckenberger-classif RS}.

\medbreak
\emph{(ii)}  \cite{A-Angiono-Heckenberger-triang}
If $V = V_1$ is a single block, then it has  a basis $(x_k)_{k\in\I_\ell}$, $\ell \ge 2$,  such that 
\begin{align}\label{equation:basis-block}
c(x_k \otimes  x_1) &= \epsilon x_1 \otimes  x_k,& c(x_k \otimes  x_j) &=(\epsilon x_j+x_{j-1}) \otimes  x_i.
\end{align}
for $k, j \in \I_\ell$, $1 < j$, where $\epsilon\in \kut$. This braided vector space is denoted by $\mathcal V(\epsilon,\ell)$. Then $\GK \toba(\mathcal V(\epsilon,\ell)) < \infty$ iff $\ell = 2$ and $\epsilon \in \{\pm 1\}$.
Let $x_{21} =  x_2x_1 + x_1 x_2$. Then
\begin{align}\label{eq:jordan-plane}
\toba(\mathcal V(1,2)) &\simeq \ku\langle x_1,x_2 \vert x_2x_1-x_1x_2+\frac{1}{2}x_1^2\rangle,
\\ \label{eq:super-jordan-plane}
\toba(\mathcal V(-1,2)) &\simeq  \ku\langle x_1, x_2 \vert x_1^2,
\, x_2x_{21}- x_{21}x_2 - x_1x_{21}\rangle
\end{align}
are the  Jordan plane  and   the super Jordan plane respectively.

\medbreak
\emph{(iii)} 
The classification of those $V$ decomposing  as  direct sum of blocks and points with $\GK \toba(V) < \infty$
appears in \cite{A-Angiono-Heckenberger-triang}; see the Introduction of \emph{loc. cit.} for details.
 
 \medbreak
 \emph{(iv)}  The classification of those $V$ with with $\GK \toba(V) < \infty$ and at least one pale block 
 was achieved in \cite{A-Angiono-Heckenberger-triang}  for $\dim V = 3$
 and in \cite{A-Angiono-Moya} for  $\dim V = 4$.  \qed
\end{theorem}
Part (i) was conjectured in \cite{A-Angiono-Heckenberger-triang} motivated by the early \cite{Rosso}. The classification in part (iv) without restrictions on $\dim V$ is work in progress.

\begin{remark}
Let  $\varGamma$ be a finite group and $V\in \yd{\ku \varGamma}$  of diagonal type
(always the case if $\varGamma$  is abelian)
such that $\GK \toba(V)  < \infty$.   
By Theorem \ref{thm:abelian-case} (i),
if $U$ is the braided vector subspace of $V$ arising from any connected component of its Dynkin diagram,
then either $\dim \toba(U)  < \infty$ or $\dim U  = 1$ with label $q=1$.
Indeed a Dynkin diagram in the list of \cite{Heckenberger-classif RS} giving rise to an infinite-dimensional Nichols algebra either is a point with label $q=1$ or has a point with label
$q\in \kut \backslash \Gb_{\infty}$,  not  realizable in $\yd{\ku \varGamma}$.
\end{remark}

\subsection{Pre-Nichols versus post-Nichols}\label{subsec:pre-post-nichols}
Here we deal with Step \ref{item:levante2}. 
Let  $(V, c)$ be a finite-dimensional braided vector space.
Recall  the notion of post-Nichols algebra, introduced  in page \pageref{text-defn:post-nichols}.
If $\Rc$ is a post-Nichols algebra of $V$, then
there is an inclusion $\toba(V)\hookrightarrow\Rc$ of graded Hopf algebras.
We complement with the following definition due to Masuoka.
A \emph{pre-Nichols algebra} of $V$ is a connected graded braided Hopf algebra $\Bc = \oplus_{n\ge 0} \Bc^n$ 
such that $\Bc^1 \simeq V$ generates $\Bc$ as algebra; thus there is an projection $\Bc \twoheadrightarrow\toba(V)$ of graded Hopf algebras.
Clearly $\toba(V)$ is the only post- and pre-Nichols algebra of $V$ up to isomorphisms.

\smallbreak
The tensor algebra $T(V)$ is naturally a pre-Nichols algebra of $V$ with the comultiplication induced by $c$.
Now the set of all pre-Nichols algebras of $V$ (up to isomorphisms) is  a poset
$\pre(V)$ with ordering given by the surjections; 
the minimal element in $\pre(V)$ is $T(V)$, and the maximal is $\Bc(V)$.

\smallbreak
Analogously, there is a graded braided Hopf algebra $T^c(V)$, which is $T(V)$ with the quantum shuffle multiplication and the natural comultiplication.
Again, the set of all post-Nichols algebras of $V$ (up to isomorphisms) is  a poset $\post(V)$  
ordered by the inclusion; now the minimal element is $\Bc(V)$ and the maximal is $T^c(V)$. If $\Omega$ is the sum of all quantum symmetrizers, then
we could picture these posets as follows:
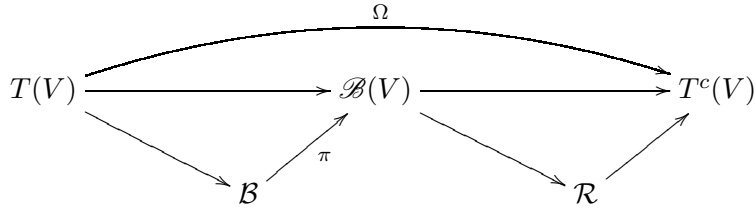
\begin{align*}
\xymatrix{T(V) \ar@/^2pc/^{\Omega}[0,6]
\ar  @{->}[rrr] \ar  @{->}[1,2] &  & & \toba(V) \ar  @{->}[1,2] \ar  @{->}[rrr] & &  & T^c(V)
\\
& &   \Bc  \ar  @{->}_{\pi}[-1,1]   &  & &  \Rc \ar  @{->}[-1,1]  &
}
\end{align*}
It is convenient to consider some relevant subposets of $\post(V)$ and $\pre(V)$.
First, if $V \in \yd{H}$; then $T(V)$, $T^c(V)$ and $\toba(V)$ are  Hopf algebras in 
$\yd{H}$; let $\post_H(V)$ and $\pre_H(V)$ be  the subposets  with elements in $\yd{H}$.
Next, $\postf(V)$ and $\pref(V)$ are the subposets whose elements have finite $\GK$.  
In this notation, Step \ref{item:levante2} amounts to determine $\postfh{H}(V) \coloneqq \postf(V) \cap \post_H(V)$ for $V \in \yd{H}$,  Question \ref{question:postnichols-finite-gkd-gral}. \phantomsection \label{text-question:postnichols-finite-gkd-gral} 
As usual we start by studying $\postf(V)$ for any braided vector space $V$.

\smallbreak
Now the graded dual $\Rc^d$ of $\Rc \in \post(V)$ is a pre-Nichols algebra of $V^*$, and vice versa. That is, the posets $\post(V)$ and $\pre(V^*)$ are anti-isomorphic. In the context of finite $\GK$,  the following holds:
 
\begin{lemma}\label{lemma:prenichols-finite-gkd-postnichols} \cite{A-Angiono-Heckenberger-Jordan}
Let $\Bc\in \pre(V)$, $\Rc = \Bc^d$.
Then  $\GK \Rc \leq \GK \Bc$. If $\Rc$ is finitely generated, then $\GK \Rc = \GK \Bc$. \qed
\end{lemma}

This leads us to consider $\pref(V)$ for any braided vector space $V$.
For this, we say as in \cite{A-Sanmarco} that $\htoba\in \pref(V)$
is \emph{eminent}  if it is a minimum in $\pref(V)$, i.e., 
\begin{itemize}[leftmargin=*]
\medbreak\item  $\GK \htoba < \infty$, and

\medbreak\item  if
$\Bc \in \pref(V)$, then there is an
epimorphism of pre-Nichols algebras $\htoba \twoheadrightarrow \Bc$.
\end{itemize}

The simplest example of vector spaces  braided by the usual flip shows that eminent pre-Nichols algebras 
are not always available, see Lemma \ref{lema:prenicholsV}. However, they do exist for many 
braided vector spaces of diagonal type with connected Dynkin diagram and they are related to the distinguished
pre-Nichols algebras $\wtoba_{\bqg}$ introduced in  \cite{Angiono-distinguished}.

\begin{theorem}\label{thm:general} \cite{A-Sanmarco,Angiono-Campagnolo-Sanmarco-sup-st, Angiono-Campagnolo-Sanmarco-supmod-ufo}
Let $\bqg$ be a braiding matrix such that $\dim\toba_{\bqg}<\infty$ and the Dynkin diagram of $\bqg$ is connected.  
\begin{itemize}[leftmargin=*]
\medbreak\item   The distinguished pre-Nichols algebra $\wtoba_{\bqg}$ is eminent except when $\bqg$ is of Cartan type
  $A_{\theta}$ or $D_{\theta}$ with $q=-1$; or 
 Cartan type $A_2$ with $q\in \Gb_3'$, or super type $\superqa{3}{q}{\{2\}}$ or $\superqa{3}{q}{\{1,2,3\}}$, 
 with $q\in \Gb_{\infty}$, or type $\g(2,3)$ with any of the diagrams
\begin{align*}
&d_1: \, \xymatrix@C=30pt{\overset{-1}{\circ} \ar@{-}[r]^{\xi^2} &\overset{\xi}{\circ}\ar@{-}[r]^{\xi}\ &\overset{-1}{\circ}}, &
&d_2: \, \xymatrix@C=30pt{\overset{-1}{\circ} \ar@{-}[r]^{\xi} &\overset{-1}{\circ}\ar@{-}[r]^{\xi}\ &\overset{-1}{\circ}}.
\end{align*}

\medbreak\item  If $\bqg$ is of type $\superqa{3}{q}{\{2\}}$ or $\superqa{3}{q}{\{1,2,3\}}$, then there is
a braided central extension of $\wtoba_{\bqg}$ by a polynomial algebra in one variable
which is eminent.

\medbreak\item  If $\bqg$ is of Cartan type  $A_2$ with $q\in \Gb_3'$, or of type $\g(2,3)$ with diagram $d_1$ or $d_2$, then
there is
a braided central extension $\htoba_{\bqg}$ of $\wtoba_{\bqg}$ by a $q$-polynomial algebra in two variables
which is eminent.  \qed
\end{itemize}
\end{theorem}

The existence of eminent pre-Nichols algebras  when $\bqg$ is of Cartan type
 $A_{\theta}$ or $D_{\theta}$ with $q=-1$ is an open problem, see Question \ref{question:postnichols-finite-gkd-diagonal}. 
\phantomsection \label{text-question:postnichols-finite-gkd-diagonal}
Eminent pre-Nichols algebras of quantum linear spaces were discussed in \cite[Theorem 1.1]{A-Sanmarco}.

\subsection{Post-Nichols algebras of vector spaces}\label{subsec:post-Nichols-vs}

Here we determine the poset $\post(V)$ of a braided vector space $(V, \tau)$, 
$\dim V < \infty$, where $\tau$ is the usual flip. 
It is well-known that $\toba(V) \simeq S(V)$, the symmetric algebra on $V$;
and that the Hopf algebra $T(V)$ is isomorphic to
the enveloping algebra $U(L(V))$,  where $L(V)$ is the free Lie algebra on $V$.
We start by the pre-Nichols algebras.

\begin{lemma}\label{lema:prenicholsV}
The pre-Nichols algebras of $V$ are of the form $U(\ngo)$, where $\ngo = \oplus_{n\in \N} \ngo^n$
is a graded Lie algebra generated by $\ngo^1 \simeq V$. Thus pre-Nichols algebras of $V$ with finite $\GK$
correspond to those $\ngo$ with $\dim \ngo < \infty$. \qed
\end{lemma}

A similar statement holds for super vector spaces.
Clearly, any finite-dimen\-sional $\ngo$ as in the Lemma is nilpotent.  
We next deal with post-Nichols algebras. The exposition requires some knowledge of algebraic groups.
We need a well-known result. 

\begin{lemma}\label{lemma:inclusion-alg-gps} \cite[13.2, 14.1]{Waterhouse}
Let $L =  \Oc(M)$ be a finitely generated commutative Hopf algebra,
$K =  \Oc(N)$ be a Hopf subalgebra of $L$ and $\pi: M \to N$ be the morphism of algebraic groups corresponding to the inclusion. Then $\pi$ is surjective. \qed
\end{lemma}

Let $\ngo$ be a finite-dimensional nilpotent Lie algebra. 
Choose a  finite-dimensional faithful representation $\ngo \to \gl(X)$ by nilpotent matrices, 
which is possible by Ado's theorem\footnote{See e.g. \texttt{https://terrytao.wordpress.com/2011/05/10/ados-theorem/}}.
Then the exponential map induces an isomorphism of varieties from $\ngo$ onto $N \coloneqq  \exp \ngo \leq \GL(X)$, which is a unipotent algebraic group with Lie $N = \ngo$.
Conversely, any unipotent algebraic group $N$ arises in this way, with $\ngo = \ln N$.

A classical result tells that a finitely generated commutative Hopf algebra with trivial coradical
is the algebra of functions on a unipotent algebraic group. 
The same holds if finite generation is replaced by finite $\GK$.

\begin{lemma}\label{lem:trivialcoradical} Let $H$ be a commutative Hopf algebra with trivial coradical. 
Assume that $\GK H = n < \infty$.
Then $H \simeq  \Oc(N)$, where $N$ is a unipotent algebraic group.
\end{lemma}

\pf Since any finite-dimensional subspace of $H$ is contained in a finite-dimensio\-nal subcoalgebra,   and hence
in an affine Hopf subalgebra, 
\begin{align*}
\GK H  =    \sup\{\GK K\vert K  \text{ is a finitely generated Hopf subalgebra}\}.
\end{align*}
Also, since $H_0 = \ku$, any finitely generated Hopf subalgebra $K$ is of the form $K \simeq  \Oc(N)$, where $N$ is a unipotent algebraic group.
Hence, in particular, $H$ is a domain.
Pick a finitely generated Hopf subalgebra  $K \simeq  \Oc(N)$ such that $\GK H = n = \GK K$. 
We claim that $H = K$. Indeed, let $x \in H$. Then $K[x]$ is contained in a finitely generated Hopf subalgebra  $L \simeq  \Oc(M)$,
and $\GK L = n$. By Lemma \ref{lemma:inclusion-alg-gps}, the map $\pi: M \to N$ is surjective; since both are unipotent, it is an isomorphism,
hence $L=K$. \epf

\begin{proposition}\label{prop:postnichols-V} Let $V$ be a finite-dimensional   vector space.
Then the finite GK-dimensional post-Nichols algebras of $V$ are of the form $ \Oc(N)$, where 
$N$ is a   unipotent algebraic group such that $\Lie N = \ngo = \oplus_{n\in \N} \ngo^n$ is a finite-dimensional graded Lie algebra generated by $\ngo^1 \simeq V^*$.
\end{proposition}

\pf Let $\Rc\in \postf(V)$ and $\toba = \Rc^d$; since $\toba$ is cocommutative, $\Rc$ is commutative. 
By Lemma \ref{lem:trivialcoradical}, $\Rc \simeq  \Oc(N)$ where $N$ is a unipotent algebraic group.
By Lemma \ref{lemma:prenichols-finite-gkd-postnichols}, $\GK \toba = \GK \Rc = \dim N$. 
By Lemma \ref{lema:prenicholsV}, $\toba = U(\ngo)$ where $\ngo = \oplus_{n\in \N} \ngo^n$ is a finite-dimensional graded Lie algebra 
generated by $\ngo^1 \simeq V^*$ and $\dim \ngo = \GK \toba $; thus $\ngo = \Lie N$. 
\epf

 \subsection{Pointed Hopf algebras over non-abelian groups}\label{subsec:Nichols-nonabelian}
 
 Let  $H$ be a pointed Hopf algebra such that 
 $\varGamma = G(H) $ is a finitely generated group. 
 If $\GK H < \infty$, then $\varGamma$ is nilpotent-by-finite. 
We are led to Question \ref{question:nichols-finite-gkd}. We have to start studying Nichols algebras with finite $\GK$ over nilpotent-by-finite groups.
First,   $V \in \yd{\ku\varGamma}$ iff $V= \oplus_{ g\in \varGamma} V_g$ is a  
$\varGamma$-graded vector space with an action of 
$\varGamma$ such that  $h \cdot V_{g} = V_{hgh^{-1}}$, $h,g \in \varGamma$.  
Then $\supp V \coloneqq \{\g \in \varGamma: V_{g} \neq 0\}$  is a subrack of $\varGamma$ with respect
to the conjugation $h\trid g = hgh^{-1}$, $h,g \in \varGamma$. 
We conclude that a braided vector space $V$ arising from a Yetter-Drinfeld module over a group
is determined by a pair $(X, \bqg)$ where $X$ is a rack and $\bqg: X \times X \to GL(U)$
is a 2-cocycle \cite{A-Grana}. Set $\toba(X, \bqg) \coloneqq \toba(V)$.
If, for simplicity, $\dim U = 1$, then the associated braiding of $V = \ku X$,
with basis $(e_x)_{x \in X}$,  is $c(e_x \otimes e_y) = \bqg_{xy} e_{x\trid x} \otimes e_x$,
$x,y \in X$. See \cite{A-Leyva,A-Carnovale-Garcia-VII,Heckenberger-Vendramin-rank>2} and their references 
 for the case $\dim \toba(X, \bqg) < \infty$. 
 
 We star our analysis of the finite $\GK$ context by the following  fact. 
 
\begin{theorem}\label{thm:infinite-conjugacy} \cite[2.6]{A-nilpotent}
Let $G$ be a finitely generated group and $V\in \yd{\ku G}$ such that $\Orb = \supp V$ is an infinite conjugacy class. Then $\GK \toba(V) \# \ku G = \infty$. \qed
\end{theorem}

A first application of Theorem \ref{thm:infinite-conjugacy} is to torsion-free nilpotent groups.

\begin{theorem} \cite[3.5]{A-nilpotent}   
Let $\varGamma$ be a finitely generated torsion-free nilpotent group
and let  $V \in \yd{\ku \varGamma}$. If  $\GK \toba(V) \# \ku \varGamma < \infty$, then $\supp V\subset Z(\varGamma)$,
thus the braided vector space $V$ `comes from the abelian case'. If  $V$  is semisimple of finite length, then its structure  can be described. \qed
\end{theorem}

By Theorem \ref{thm:infinite-conjugacy} we are reduced
to Question \ref{question:racks-finite-gkd}. \phantomsection \label{text-question:nichols-finite-gkd}
We next discuss the role of racks of type C.
 Translating  \cite[Theorem 2.1]{Heckenberger-Vendramin-rank2} to the language of racks, the notion of
 finite rack 
 of type C was introduced in \cite{A-Carnovale-Garcia-III}. It comes with the following  criterium.
 
 \begin{theorem} 
\cite{A-Carnovale-Garcia-III}  A finite rack $\Orb$ of type C satisfies
 $\dim\toba(\Orb, \bqg) = \infty$ for every 2-cocycle $\bqg$. \qed
 \end{theorem}

 By analogy with the abelian case  \cite{A-Angiono-Heckenberger-triang,Angiono-Garcia2}, 
 Conjecture \ref{conjecture:typeC-X-finite-GK} proposes extending the preceding statement to $\GK$.
 \phantomsection \label{text-conjecture:typeC-X-finite-GK} Indeed, one wonders whether 
 \cite[Theorem 2.1]{Heckenberger-Vendramin-rank2}, proved using the Weyl groupoid, holds for finite $\GK$.
The validity of  Conjecture \ref{conjecture:typeC-X-finite-GK} would have a wide impact, for instance:

\begin{itemize}[leftmargin=*] 
\medbreak\item[$\circ$]  \cite[2.9]{A-nilpotent}
If $\varGamma$ is a finitely generated  nilpotent group 
with torsion subgroup $T$ and $1 < \vert T\vert$ is  odd, then
a finite conjugacy class of $\varGamma$ is either abelian or  of type C. 
\end{itemize}

There is a long list of simple racks, mostly conjugacy classes of simple groups, that are of type C, 
see \cite{A-Leyva,A-Carnovale-Garcia-VII} and their references. But there are racks with rather uncomplicated 
structure that are not of type C and give rise  to difficult questions.
Let $X$ be a finite  indecomposable rack with at least 2 elements
and consider the  cocycle $\pmb{\lambda}: X\times X\to \kut$ that takes the constant value $\lambda$.
If $\lambda \notin \Gb_{\infty}$, then $\GK\toba(X, \mathbf{1}) = \infty$ by \cite[Example 1.10.10]{Heckenberger-Schneider-book} but for $\lambda \in \Gb_{\infty}$ there are no clues; 
the important case $\lambda = 1$ is recorded as
Question \ref{question:nichols-cocycle-1}. \phantomsection\label{text-question:nichols-cocycle-1}

\subsection{Noetherian  pointed Hopf algebras}\label{subsec:pointed-noetherian}
Since the classical Questions  \ref{question:enveloping-Noetherian} and 
\ref{question:group-algebra-Noetherian}
are open, cf. also the observations in Remark \ref{rem:text-noetherian}, 
it is hard to imagine any approach to the classification of Noetherian  pointed Hopf algebras.
But  some comments are in order.
Let $H$ be a pointed Hopf algebra with $\varGamma \coloneqq G(H)$  finitely generated, 
so that $\gr H \simeq \Rc \# \ku \varGamma$ where $\Rc$ is the diagram of $H$.

\begin{itemize}[leftmargin=*]
\medbreak\item[$\circ$]  If $\varGamma$ is polycyclic-by-finite and $R$ a Noetherian Hopf algebra in $\yd{\ku \varGamma}$,
then $R \# \ku \varGamma$ is Noetherian, see \cite[1.5.12]{Mcconell-Robson}.

\medbreak\item[$\circ$]  If $\GK H < \infty$, then $\varGamma$ is nilpotent-by-finite, hence polycyclic-by-finite. Thus,
if $\Rc$ is Noetherian, then $\gr H$ and a fortiori $H$ are Noetherian, by the previous claim and Lemma \ref{lemma:filtration-noeth}.

\medbreak\item[$\circ$]  The Nichols algebras described in \cite{A-Angiono-Heckenberger-triang,A-Angiono-Moya}  have a convex PBW-basis, thus are Noetherian; so are their bosonizations
with polycyclic-by-finite group algebras. 

\medbreak\item[$\circ$]  (Heckenberger, private communication). 
The coideal left ideals of a pointed brai\-ded Hopf algebra are
in bijection with the left coideal subalgebras. Hence, in order to prove
that a Nichols algebra is not noetherian, it suffices to construct
an infinite strictly increasing chain of left coideal subalgebras. For
Nichols algebras of diagonal type,   
such an increasing chain exists iff the root system is not arithmetic.
\end{itemize}

\subsection{(Pointed) Hopf algebras with low $\GK$}\label{subsec:pointed-lowGk}
The problem of classifying Hopf algebras with $\GK =1$ or $2$ was addressed by several authors
under suitable restrictions. See \cite{Brown-Zhang-survey} for a recent report. We highlight:
 \begin{itemize}[leftmargin=*]
\medbreak\item[$\circ$]  The classification of affine prime regular Hopf algebras with $\GK =1$ was  completed in \cite{Wu-Liu-Ding}; there are non-pointed examples, answering  a question in \cite{Brown-Zhang}.

\medbreak\item[$\circ$]  The classification was  extended in \cite{Liu} relaxing the hypothesis of regularity.
 \end{itemize}
 
\section{Beyond pointed Hopf algebras}\label{section:more}

We discuss two applications and two generalizations of the lifting method.
Throughout $H$ is a Hopf algebra. 

\subsection{Nichols algebras over algebraic groups}\label{subsec:algebraic-groups}

The next instance of the lifting method  is to consider \emph{commutative} 
cosemisimple Hopf algebras  in Step \ref{item:levante0}.  
We begin by looking at the finite-dimensional simple Yetter-Drinfeld modules
over an affine commutative Hopf algebra $K$. 
\emph{Copointed} Hopf algebras, i.e., with coradical $\simeq \ku^{G}$ where $G$ is
a finite group, were studied e.g. in \cite{Garcia-Vay}, so we assume that $0 < \GK K$. 
Let $\Gb \coloneqq \Alg(K, \ku)$ be the corresponding affine algebraic group so that $K \simeq \Oc(\Gb)$, cf. the notation in page \pageref{text-def:simple-algebraic}.
Before starting recall the induced representation in the setting of algebraic groups.
Given an algebraic subgroup $\Hb \leq \Gb$, the functor $\Ind_{\Hb}^{\Gb}$ assigns to $M \in \rep \Hb$ the induced (rational) module
\begin{align*}
\Ind_{\Hb}^{\Gb} M \coloneqq (M \otimes \Oc(\Gb))^{\Hb} = \{\varphi \in \Hom_{\text{var}}(\Gb, M): \varphi(gh) = h^{-1}\varphi(g)\, \forall g\in \Gb, h \in \Hb\}
\end{align*}
where $\Gb$ acts on $\Ind_{\Hb}^{\Gb} M$ by $(\gamma \cdot \varphi)(g) = \varphi(\gamma^{-1}g)$, 
$\gamma, g \in \Gb$, $\varphi \in \Ind_{\Hb}^{\Gb} M$; see \cite[3.3]{Jantzen}.

\medbreak
Let $\gamma  \in FC(\Gb)$ and  $d \coloneqq |\Orb_{\gamma}| < \infty$. Fix  a 
family  $(g_{\kappa})_{\kappa\in \Orb_\gamma}$ in $\Gb$ such that 
$\kappa = g_{\kappa} \trid \gamma$, $\kappa\in \Orb_\gamma$.
Thus $\displaystyle\Gb \overset{\star}{=} \varcoprod_{\kappa\in \Orb_\gamma} g_{\kappa}\Cg{\gamma}$ where $g_{\kappa}\Cg{\gamma}$ is closed and open.

Let  $U \in\rep \Cg{\gamma}$.   By $\star$, the induced  module $\Mt(\gamma,U)  \coloneqq \Ind_{\Cg{\gamma}}^{\Gb} U $ has a grading of vector spaces $\Mt(\gamma,U) = \oplus_{\kappa\in \Orb_\gamma} \Mt(\gamma,U)^{\kappa}$, where
\begin{align}\label{eq:grading-induced}
\Mt(\gamma,U)^{\kappa} &= \{\varphi \in \Ind_{\Cg{\gamma}}^{\Gb} U: 
\varphi(g_{\eta}h) = 0\, 
\forall \eta\neq \kappa, h \in \Cg{\gamma}\}. 
\end{align}
Clearly,
$\dim	\Mt(\gamma,U)  = \vert \Orb_{\gamma} \vert \, \dim U$.

\begin{lemma} \begin{enumerate}[leftmargin=*,label=\textrm{(\alph*)}]
\medbreak\item \label{item:YD-alg-gps-a}
Let $V \in \yd{K}$, $\dim V < \infty$. 
Then $\supp V \subseteq FC(\Gb)$ is stable by conjugation of $\Gb$
and $\soc V  = \oplus_{\kappa \in \Gb} V^{\kappa}$ is a Yetter-Drinfeld 
submodule\footnote{Here $\soc V$ is the socle of $V$  in $\lmod{K}$.} 
of $V$.
If $V$ is simple, then $V = \soc V$
 and $\supp V$ is  a conjugacy class  of $\Gb$.

\medbreak\item \label{item:YD-alg-gps-b} Let $V \in \rep \Gb$, $\dim V < \infty$. Assume that
$V$ has a grading $V = \oplus_{\kappa\in \Gb} V^{\kappa}$
so that $V$ is a semisimple $\Oc(\Gb)$-module. Then $V \in \yd{K}$ if and only if
\begin{align}\label{eq:alg-gps-yd-ss}
\eta \cdot V^{\kappa} &= V^{ \eta \trid \kappa}& \text{for all  } \eta,\kappa&\in \Gb.
\end{align}

\medbreak\item \label{item:YD-alg-gps-c} Let $\gamma  \in FC(\Gb)$ and $U \in\rep \Cg{\gamma}$.   The induced  module $\Mt(\gamma,U)  \coloneqq \Ind_{\Cg{\gamma}}^{\Gb} U $ with the action of $\Oc(\Gb)$ given by the grading \eqref{eq:grading-induced}  belongs to $\yd{H}$.

\medbreak\item \label{item:YD-alg-gps-d} If $U \in \irr \Cg{\kappa}$, then 
$\Mt(\kappa,U)$ is a simple finite-dimensional object in $\yd{H}$ and any simple
finite-dimensional object in $\yd{K}$ arises in this way.
\end{enumerate}
\end{lemma}

\pf    \ref{item:YD-alg-gps-a}, \ref{item:YD-alg-gps-b}: By \eqref{eq:action-rational}, 
$V\in \rep \Gb$; since $V \in \rep \Oc(\Gb)$,  
$V = \oplus_{\gamma\in \Gb} V^{(\gamma)}$, cf. \eqref{eq:decomp-commutative-dix}. 
In terms of the action of $\Gb$, the compatibility  of Yetter-Drinfeld modules  reads
\begin{align*}
& &\gamma \cdot (f\cdot v)  
&= \langle \gamma^{-1}, f\_1\rangle \langle \gamma , f\_3\rangle \,
f\_2\cdot (\gamma \cdot v),  
\end{align*}
for all $v\in V$, $\gamma \in \Gb$, $f\in \Oc(\Gb)$, which is equivalent to 
\begin{align*}
f\cdot v  &= \langle \gamma^{-1}, f\_1\rangle \langle \gamma , f\_3\rangle \,
\gamma^{-1} \cdot \big(f\_2\cdot (\gamma \cdot v)\big)
\end{align*}
or setting $w \coloneqq \gamma \cdot v$,  to $f\cdot \big( \gamma^{-1} \cdot w \big) = \langle \gamma^{-1}, f\_1\rangle \langle \gamma , f\_3\rangle \,
\gamma^{-1} \cdot  (f\_2\cdot w)$, i.e., to
\begin{align}\label{eq:alg-gps-yd-ss2}
f\cdot \big( \eta\cdot w \big) &= \langle \eta, f\_1\rangle \langle \eta^{-1} , f\_3\rangle \,
\eta \cdot  (f\_2\cdot w),&
\forall w\in V, \eta &\in \Gb, f\in \Oc(\Gb).
\end{align}
Now, if $\kappa\in \Gb$ and $w\in V^{\kappa}$,  then \eqref{eq:alg-gps-yd-ss2}
is equivalent to \eqref{eq:alg-gps-yd-ss}, implying \ref{item:YD-alg-gps-b}.
Thus, for any $\gamma, \kappa\in \Gb$,  $\gamma \cdot V^{\kappa} = V^{\gamma  \trid \kappa}$
and $\gamma \trid  \supp V = \supp V$. 
Hence $\soc V  = \oplus_{\kappa \in \Gb} V^{\kappa}$ is a Yetter-Drinfeld submodule of $V$
and  $\supp V \subseteq FC(\Gb)$.
For a conjugacy class $\Orb \subseteq FC(\Gb)$, $V^{\Orb} \coloneqq \oplus_{\gamma\in \Orb} V^{\gamma}$,  is a Yetter-Drinfeld submodule of 
$\soc V = \oplus_{\Orb \subseteq FC(\Gb)} V^{\Orb}$; if $\soc V$ is indecomposable, then 
$\soc V = V^{\Orb}$ for some   $\Orb \subseteq FC(\Gb)$. This shows \ref{item:YD-alg-gps-a}.

\ref{item:YD-alg-gps-c} Given $g \in \Gb$ and $\eta, \kappa \in \Orb_{\gamma}$, 
we have 
\begin{align*}
 g^{-1}g_{\eta} \trid \gamma  = g^{-1} \trid \eta  &
 \implies g^{-1}g_{\eta} \in g_{g^{-1} \trid \eta}  \Cg{\gamma}.
\end{align*}
 Thus, if $h \in \Cg{\gamma}$ and $\varphi \in \Mt(\gamma,U)^{\kappa}$, then
 for some $h' \in  \Cg{\gamma}$ we have
\begin{align*}
 g\cdot \varphi(g_{\eta}h)  &= \varphi(g^{-1}g_{\eta}h) 
 = \varphi(g_{g^{-1} \trid \eta}h') \neq 0 \iff  g^{-1} \trid \eta = \kappa
 \iff  \eta = g \trid\kappa.
\end{align*}
That is, $g\cdot \varphi \in \Mt(\gamma,U)^{g \trid\kappa}$ and \eqref{eq:alg-gps-yd-ss} holds.

\ref{item:YD-alg-gps-d} Let $V \in \yd{K}$, $\dim V < \infty$. 
Fix  $\gamma  \in \supp V$; then $V \in \rep \Cg{\gamma}$ (the restriction of a rational module to a closed subgroup is rational) and
 $V^{\gamma} \in \rep \Cg{\gamma}$ (submodules of  rational modules are rational). 
 Now $\sum_{g \in \Gb} g \cdot U$ is a Yetter-Drinfeld submodule of $V$ isomorphic to $M(\kappa,U)$ for any $\Cg{\gamma}$-submodule $U$ of $V^{\gamma}$.
This implies the claim.  \epf

If $\kappa  \in FC(\Gb)$, then $\Cg{\kappa}$ contains the connected component $\Gb^0$ of $\Gb$. 
\emph{Assume now that $\Gb$ is connected}, so that $FC(\Gb) = Z(\Gb)$. 
Then the simple objects in $\yd{K}$ are of the form $U^{\kappa} = U\in\irr{\Gb}$,
homogeneous of degree $\kappa  \in Z(\Gb)$. 
By the Schur Lemma, $\kappa$ acts by a scalar $q^{-1} \in \kut$ on $U$. 
Then the braiding of $U^{\kappa}$ is given by
\begin{align*}
c(v\otimes w) &= v\_{-1}\cdot w \otimes v\_0 =  \langle \kappa, v\_{-1} \rangle  w \otimes v\_0 
= w \otimes \kappa^{-1} \cdot v  = q  w\otimes v, & v,w&\in U.
\end{align*}
Thus $U^{\kappa}$ is a braided vector space of diagonal type with matrix $q_{ij} = q$ for all $i,j\in \I_{ n}$, $n = \dim U$.  
By Theorem \ref{thm:abelian-case} (i),  $\GK\toba(U^{\kappa}) < \infty$  iff
$q$ and $n$ are as in Table \ref{tab:nichols-q-n}.

\begin{table}[ht]
	\caption{$\GK\toba(U^{\kappa}) < \infty$}\label{tab:nichols-q-n}
	\begin{tabular}{| c | c | c |c |c | }
\hline
$q$ & $n$ & $\toba(U^{\kappa})$ & $\GK$ & $\dim$
\\  \hline
$\in \Gb_{N}'$, $N > 1$ & $1$ & $\ku[X]/(X^N)$ & & $N$
\\ \hline 
$\in \Gb_{3}'$ & $2$ & type $A_2$ &  & $27$
\\ \hline
$-1$ & arbitrary & $\Lambda (U)$ & & $2^n$
\\ \hline
$1$ & arbitrary & $S (U)$ & $n$ &
\\ \hline
$\in \ku \backslash\Gb_{\infty}$ & $n=1$ & $\ku[X]$ & 1 & 
\\ \hline 
\end{tabular}
\end{table}

\medbreak 
Assume next that $V \in \yd{K}$ is semisimple, thus
$V = \oplus_{j \in \I_\theta}U_j^{\kappa_j}$, where $\kappa_j  \in Z(\Gb)$, 
$U_j \in\irr{\Gb}$,  $n_j = \dim U_j$.  By the Schur Lemma, $\kappa_j$ acts by 
$q_{ij}^{-1}\in \kut$ on $U_i$; hence $V$ is of diagonal type and 
its braiding matrix is obtained gluing $(i,j)$-blocks of the form $q_{ij} \id \in \ku^{n_i\times n_j}$.  In each case 
it is possible to decide whether  $\GK \toba(V) < \infty$.

\begin{remark}
There might be $V \in \yd{K}$ non-semisimple with $\GK \toba(V) < \infty$.
For example, if 
$\Gb \simeq \Tb$, a torus of dimension $d$, then $\Oc(\Tb) \simeq \ku \Z^d$, 
treated in Theorem \ref{thm:abelian-case}. 
\end{remark}

\smallbreak
We now focus on the initial goal of this Subsection, namely affine Nichols algebras 
over $\Oc(\Gb)$ cosemisimple, or equivalently $\Gb$ reductive. 
We start by an example.

\begin{example}  \label{exa:nichols-sl2}
	Let $\Gb = SL(2, \ku)$ and $L(n)$ the simple $\Gb$-module of dimension $n+1$. Then  
$Z(\Gb) = \{\pm \id\}$ acts trivially on $L(n)$ iff $n$ is even. 
Let $V = \oplus_{j \in \I_\theta}L(n_j)^{\kappa_j}$ and set 
\begin{align*}
\I_- &= \{t \in \I_{\theta}: \kappa_t = -\id, \ n_t  \text{ is odd} \},&
\I_+ &= \{t \in \I_{\theta}: \kappa_t = \id, \ n_t  \text{ is even} \}.
\end{align*}
Then  $\GK \toba(V) < \infty$ iff either of the following holds:
\begin{enumerate}[leftmargin=*,label=\textrm{(\alph*)}]
\medbreak\item \label{item:exa-nichols-sl2-a} $\I_{\theta} = \I_- \coprod \I_+$. Set $V_{\pm} = \oplus_{j \in \I_{\pm}}L(n_j)^{\kappa_j}$.  In this case,  
$\toba(V) \simeq \Lambda(V_-) \otimes S(V_+)$ and $\GK \toba(V) \# \Oc(SL(2, \ku)) = \dim V_+  + 3$.

\medbreak\item \label{item:exa-nichols-sl2-b} There exists $j \in \I_{\theta}$ such that $\kappa_j = -\id$ and $n_j$ is even;
then $n_i$ is even for all $i \in \I_{\theta}$. In this case, 
$\toba(V) \simeq S(V)$  and $\GK \toba(V) \# \Oc(SL(2, \ku)) = \dim V + 3$.

\medbreak\item \label{item:exa-nichols-sl2-c} There exists $j \in \I_{\theta}$ such that $\kappa_j = \id$ and $n_j$ is odd;
then $\kappa_i = \id$ for all $i \in \I_{\theta}$. In this case, 
$\toba(V) \simeq S(V)$  and $\GK \toba(V) \# \Oc(SL(2, \ku)) = \dim V + 3$.
\end{enumerate}
\pf When $\I_{\theta} = \I_- \coprod \I_+$, this is clear. Otherwise, 
there exists $j \in \I_{\theta}$ either as in \ref{item:exa-nichols-sl2-b} or as in \ref{item:exa-nichols-sl2-c}. Now if $W$ is a braided vector space of diagonal type with Dynkin diagram $\xymatrix{\overset{1}{\circ} \ar  @{-}[r]^{ r}  &  \overset{s}{\circ}}$ where 
$r \neq 1$, then $\GK W = \infty$, cf. \cite{A-Angiono-Heckenberger-triang,Rosso};
this implies the claims. 
\epf
For instance, if $\I = \I_-$, then $\toba(V) \simeq \Lambda (V)$ and 
$\toba(V) \#  \Oc(\Gb)$ is co-Frobenius, but not cosemisimple.
Similarly, the only co-Frobenius Hopf algebra whose coradical is $\Oc(PSL(2, \ku))$ 
where  $PSL(2, \ku) = \Gb/ Z(\Gb)$, is $\Oc(PSL(2, \ku))$ itself.
\end{example}

Example \ref{exa:nichols-sl2} illustrates the general case. 
Let $\Gb$ be a simple algebraic group with lattice of weights $\varLambda$, see page \pageref{text-notation:simple-algebraic}.
Then $\Hom(Z(\Gb), \ku^{\times}) \simeq \varLambda/Q$;  
given $L(\lambda) \in \irr \Gb$, $\lambda \in \varLambda^+$, and $\kappa\in Z(\Gb)$, 
the pairing between $\kappa$ and $\lambda \mod Q$ is the scalar 
$\chi_{\lambda}(\kappa)$ by which $\kappa$ acts on $L(\lambda)$.  
The groups $P/Q$ are well-known, see Table \ref{tab:nichols-simplegroup}.
We see that $\toba (L(\lambda)^{\kappa})$ is either the symmetric or 
the exterior algebra; 
indeed, various cases in Table \ref{tab:nichols-q-n} do not arise: 
\begin{itemize}[leftmargin=*] 
\medbreak\item[$\circ$]  If $n = 1$, then $L(\lambda)$ is trivial, hence $\chi_{\lambda}(\kappa) = 1$.

\medbreak\item[$\circ$]  If $n = 2$, then $\Gb\simeq SL(2, \ku)$ and $\chi_{\lambda}(\kappa) \in \{\pm 1\}$. 
\end{itemize}
Thus the Nichols algebras of the semisimple $V\in \yd{\Oc(\Gb)}$ 
belong to the classical or super worlds.
The previous discussion suggests that the determination of the affine Nichols algebras in
$\yd{\Oc(\Gb)}$ where now $\Gb$ is semisimple 
(or even reductive) is reachable, once
part (iv) of Theorem \ref{thm:abelian-case} is completed without  dimension constraints.

\begin{table}[t]
	\caption{Centers of simple algebraic groups}\label{tab:nichols-simplegroup}
	\begin{tabular}{| c | c | c |c |c |c | }
\hline
Type & $P/Q$ & Type & $P/Q$ & Type & $P/Q$
\\  \hline
$A_{\ell}$, $\ell \geq 1$ & $C_{\ell  +1}$ & $B_{\ell}$, 
$C_{\ell}$, $\ell \geq 2$, $E_{7}$ & $C_ 2$ & $D_{\ell}$, $\ell \geq 4$ even 
& $C_ 2 \times C_2$
\\ \hline 
$D_{\ell}$, $\ell > 3$ odd & $C_ 4$ &
$E_{6}$ & $C_ 3$&  $E_{8}$, $F_4$, $G_2$ & trivial
\\ \hline 
\end{tabular}
\end{table}

\subsection{Nichols algebras over quantum groups}\label{subsec:quantum-groups}

Another possibility in Step \ref{item:levante0} of the lifting method  is to consider the 
cosemisimple Hopf algebras $\Oc_{q}(\Gb)$ where  $\Gb$ is a simple algebraic group
and $q\in \kut \backslash \Gb_{\infty}$. 
To our knowledge their Nichols algebras were not completely described but some information
is available. Indeed, these Hopf algebras are co-quasitriangular,
thus there is an embeding of braided tensor categories 
$\lcomod{\Oc_{q}(\Gb)} \hookrightarrow \yd{\Oc_{q}(\Gb)}$. 
Now by its very definition, $\lcomod{\Oc_{q}(\Gb)}$ is equivalent to the semisimple category
of type 1 representations of $U_q(\g)$; but $\rep U_q(\g)$ is closely related to
$\yd{U_q(\bg)}$ where $\bg$ is the  Borel subalgebra of $\g$.
The Nichols algebras in $\yd{U_q(\bg)}$ can be determined by a rondabout argument sometimes
called the splitting technique that we outline next. \phantomsection\label{text:splitting-technique}
Let $H$ be a Hopf algebra, $\Vc, \Wc\in \yd{H}$ and $\Uc = \Vc \oplus \Wc$.
Then there are Hopf algebra  maps 
$\xymatrix{\toba(\Uc)\# H \ar@/^0.4pc/^{\pi}[0,1]& \toba(\Vc)\# H \ar@/^0.4pc/_{\iota}[0,-1]}$ 
such that $\pi\iota 
= \id_{\toba(\Vc)\# H}$. Therefore 
\begin{align*}
\toba(\Uc)\# H & \overset{\star}{\simeq} \K\# \toba(\Vc)\# H;&
&\text{here } &\K =  \toba(\Uc)\# H ^{\co \pi} &\simeq \toba (Z_\Wc)
\end{align*}
where	$Z_\Wc \coloneqq \ad_c\toba (\Vc) (\Wc)$.
See e.g. \cite{Heckenberger-Schneider-book}. 
The isomorphism $\star$ is used in the definition of the Weyl groupoid,  \cite{Heckenberger-Yamane,A-Heckenberger-Schneider,Heckenberger-Schneider-book};
to get information on $\toba(\Uc)$ from knowledge of $\toba(Z_\Wc)$  and $\toba(\Vc)$ 
\cite{Rosso,A-Angiono-Heckenberger-triang};
and to get information on $\toba (Z_\Wc)$ from knowledge of  $\toba (\Uc)$ and $\toba(\Vc)$ \cite{A-Angiono-basic,Rosso}. Clearly, $	\dim  \toba(\Uc) = \dim \toba (Z_\Wc)\dim  \toba(\Vc)$ while 
\begin{align*} 
 \GK  \toba(\Uc) &\leq \GK\toba (Z_\Wc) + \GK  \toba(\Vc),
\end{align*}
having equality  when $\toba(\Vc)$ has a convex PBW-basis. 
See Question \ref{question:braided-tensor-product}. \phantomsection\label{text-question:braided-tensor-product}
Precise information on $\toba (Z_\Wc)$ is available when $\Uc$ is of diagonal type  \cite{A-Angiono-basic}; this would lead to the determination of the finite GK-dimensional Nichols algebras in $\yd{U_q(\bg)}$.

\subsection{Compact  Hopf algebras}\label{subsec:compact}
In this subsection, the base field is $\Cb$. Recall that a $\ast$-algebra is an associative algebra $A$
provided with a sesquilinear involution $x \mapsto x^{\ast}$ such that
$(xy)^{\ast} = y^{\ast}x^{\ast}$, for all $x,y \in A$. If $A$ and $B$ are $\ast$-algebras, then
$A \otimes B$ is a $\ast$-algebra with 
$(x \otimes y)^{\ast} = x^{\ast} \otimes y^{\ast}$  for all $x\in A$, $y\in B$. Evidently, 
$\Cb$ is a $\ast$-algebra with $\lambda^{\ast} = \overline{\lambda}$ (complex conjugation).

\begin{itemize}[leftmargin=*]
\medbreak\item    A $\ast$-Hopf algebra is a Hopf algebra $H$ which is also a $\ast$-algebra, such that
the comultiplication and the counit are morphisms of $\ast$-algebras.
\end{itemize}
Let  $H$ be a $\ast$-Hopf algebra. Then $H^{\ast}$ is a $\ast$-algebra via
$\langle f^{\ast}, x \rangle = \overline{\langle f, (\Ss x)^{\ast} \rangle}$, for any  $x\in H$, $f\in H^*$.
It follows from the uniqueness of the antipode that  $(\Ss x)^{\ast} = \Ss^{-1}(x^{\ast})$.  

\begin{itemize}[leftmargin=*] 
\medbreak\item[$\circ$]  The following data are equivalent: (a) a right integral $\smallint: H \to \ku$,
(b) a sesquilinear form $\Lbrack \, \vert \, \Rbrack: H \times H \to \Cb$ satisfying  
\begin{align}\label{eq:killing3}
\Lbrack xy\vert z\Rbrack  &= \Lbrack y\vert x^{\ast}z\Rbrack, & x,y,z &\in H,
\\
\label{eq:killing4}
\Lbrack f \rightharpoonup x\vert y \Rbrack 
&= \Lbrack  x\vert  f^{\ast} \rightharpoonup y \Rbrack, & x,y &\in H, f\in H^*.
\end{align}
\end{itemize}
Explicitly, the correspondence is given by
$\Lbrack x\vert y \Rbrack  = (y^{\ast} \vert x) 
= \langle \smallint, y^{\ast} x \rangle$, $\langle \smallint, x\rangle = (x \vert 1) = \Lbrack 1\vert x^{\ast} \Rbrack$ for any $x,y\in H$,  cf. \eqref{eq:killing1} and \eqref{eq:killing2}. 

\begin{definition} \cite{Woronowicz-CQG}
A \emph{compact Hopf algebra} $H$ is a $\ast$-Hopf algebra 
such that any  finite-dimensional $H^*$-comodule carries an invariant inner product.
\end{definition} 

Actually this is the definition of "algebra of regular functions" on a compact quantum group (CQG)
\cite{Woronowicz-CQG}, but
it is known that the completion of a compact Hopf algebra 
with respect to a suitable norm  is a CQG  in the sense of \cite{Woronowicz-CQG}.
Clearly a compact Hopf algebra is cosemisimple. 

\begin{itemize}[leftmargin=*]
\medbreak\item[$\circ$]  Let  $H$ be a cosemisimple $\ast$-Hopf algebra and let $\smallint \in \Ic_r(H)$ with  $\langle \smallint, 1\rangle = 1$.
The following are equivalent: 
\begin{enumerate} [leftmargin=*,label=\textrm{(\alph*)}]
\medbreak\item  $H$ is a compact Hopf algebra.

\medbreak\item  The sesquilinear form $\Lbrack \, \vert \, \Rbrack$ is an inner product, i.e., it is positive definite.  
\end{enumerate}

In this case, we  say that $\ast$ is a compact involution.  

\medbreak\item[$\circ$]  \cite{A-Czechoslovak} A cosemisimple Hopf algebra has at most one compact involution modulo automorphisms.
\end{itemize}

The problem of classifying compact Hopf algebras under finiteness conditions,  
Question \ref{question:CQG}, arises naturally. 
\phantomsection\label{text-question:CQG} The  approach of classifying  cosemisimple Hopf algebras 
with finiteness conditions, Question \ref{question:cofrob-coss}, 
and then determine which one supports a compact involution seems too coarse\footnote{There are cosemisimple Hopf algebras that do not admit a compact involution.}.
See \cite{Neshveyev-Yamashita} for a strategy using analytical tools. 
A first approximation goes through the following notion.

\begin{definition}\label{def:CQG} \cite{Woronowicz-newquantum}
A cosemisimple Hopf algebra  $H$ is a \emph{quantum group of Lie type}
if there exists  a simple algebraic group  $\Gb$ with lattice $\varLambda$, as in page \pageref{text-def:simple-algebraic},   having the same fusion rules and dimensions of representations
as $H$. That is, $\irr \rcomod{H}$ is parametrized by  $\rep \Gb \simeq \varLambda^+$; 
if  $\Lc(\lambda)$ is the simple $H$-comodule corresponding to $\lambda \in \varLambda^+$, then  
$\dim \Lc(\lambda) = d_{\lambda}$ and, recalling \eqref{eq:simple-modules}, 
\begin{align*}
\Lc(\lambda) \otimes \Lc(\mu) &\simeq \oplus_{\nu \in \varLambda^+} \Lc(\nu)^{\fusion{\lambda}{\mu}{\nu}},& \lambda, \mu &\in \varLambda^+.
\end{align*}
We shall also say that $H$ is a quantum $\Gb$. We set $\nu \prec \lambda  +\mu + \dots + \xi$ 
to express that
$L(\nu)$ appears in $L(\lambda) \otimes L(\mu) \dots \otimes L(\xi)$, 
In this case $\nu \leq \lambda + \mu + \dots + \xi $, i.e., $\lambda + \mu + \dots + \xi - \nu \in Q^{+}$. 
For instance $\nu \prec \lambda +\mu$ iff $\fusion{\lambda}{\mu}{\nu} \neq 0$.
Given $\lambda \in \varLambda^+$, let $q_0(-\lambda) \in \varLambda^+$ be such that
$L(\lambda)^* \simeq L(q_0(-\lambda))$. Then $0 \prec \lambda + q_0(-\lambda)$, thus
\begin{align}\label{eq:quantumG-dual}
\lambda + q_0(-\lambda) \in Q^{+}.
\end{align} 
\end{definition} 

Let $C_{V}\hookrightarrow \Oc(\Gb)$, respectively $D_{W}\hookrightarrow H$, 
be the subcoalgebra of the matrix coefficients of $V \in \rcomod{\Oc(\Gb)}$, respectively  
$W \in \rcomod{H}$.

\begin{lemma}\label{lemma:CQG} Let  $\Gb$ be a simple algebraic group  with lattice $\varLambda$ and 
let $H$ be a quantum $\Gb$. Keep the notation above.
\begin{enumerate} [leftmargin=*,label=\textrm{(\roman*)}]
\medbreak\item \label{item:CQG-GKdim} \cite[Theorem 2.9]{Chirvasitu-Walton-Wang}
 $\GK H = \dim \Gb$.

\medbreak\item \label{item:CQG-lattice} Given a lattice $\varTheta$, $Q \leq \varTheta  \leq \varLambda$,
$H_{\varTheta} \coloneqq \oplus_{\lambda \in \varTheta^+} D_{\Lc(\lambda)}$
is a normal Hopf subalgebra of $H$ and a quantum $\Gb_{\varTheta} \coloneqq \Gb/ (\varTheta/Q)$.
\end{enumerate}
\end{lemma} 

\pf 
\ref{item:CQG-GKdim} We have $C_{L(\lambda)}C_{L(\mu)} 
= C_{L(\lambda)\otimes L(\mu)} =
\oplus_{\nu  \prec  \lambda + \mu} C_{L(\nu)}$, hence 
\begin{align}\label{eq:quantumG-tensor}
D_{\Lc(\lambda)}D_{\Lc(\mu)} = D_{\Lc(\lambda)\otimes \Lc(\mu)} =
\oplus_{\nu  \prec  \lambda + \mu} D_{\Lc(\nu)}.
\end{align}
Therefore $\dim C_{L(\lambda)}C_{L(\mu)} = \dim D_{\Lc(\lambda)}D_{\Lc(\mu)}$.
Take now $V= \oplus_{1\le j \le \rk  \g} L(\epsilon_j) $, where the $\epsilon_j$'s are the
fundamental weights, and $W= \oplus_{1\le j \le \rk  \g} \Lc(\epsilon_j) $. By the preceding claim,
$\dim V^n = \dim W^{n}$ for any $n \in \N$, hence 
$\GK H =\GK \Oc(\Gb) = \dim \Gb$.
\ref{item:CQG-lattice} $H_{\varTheta}$ is a Hopf subalgebra  by \eqref{eq:quantumG-dual} 
and \eqref{eq:quantumG-tensor}. If $x\in D_{\Lc(\lambda)}$ and $y \in D_{\Lc(\mu)}$, then 
\begin{align*}
\ad(x) (y) = x\_{1}y \Ss(x\_{2}) \in D_{\Lc(\lambda)} D_{\Lc(\mu)} D_{\Lc(q_0(-\lambda))}
= \oplus_{\nu  \prec  \lambda + \mu + q_0(-\lambda)} D_{\Lc(\nu)}.
\end{align*}
If $\mu \in \varTheta^+$, then $\lambda + \mu + q_0(-\lambda) \in \varTheta^+$ by 
\eqref{eq:quantumG-dual}, hence $H_{\varTheta}$ is stable by the adjoint action of $H$.
\epf

Now the Question \ref{question:cofrob-coss} could be
split into Questions \ref{question:CQG-Lie} and \ref{question:CQG-Lie-exhaustion}. \phantomsection\label{text-question:CQG-Lie-exhaustion} So far it is known:

\begin{itemize}[leftmargin=*] 
\item[$\circ$]  \cite{Woronowicz-newquantum} The list of the quantum $SL(2, \Cb)$'s consists of
$\Oc\left(SL(2, \Cb)\right)$, $\Oc_{q}\left(SL(2, \Cb)\right)$ where $q \notin \Gb_{\infty}$
and a Jordanian $\Oc_{J}\left(SL(2, \Cb)\right)$ described in \cite{Zakrzewski,Demidov}.

\medbreak\item[$\circ$]  The classification of the quantum $SL(3, \Cb)$'s  is given in \cite{Ohn-JA,Ohn-missing}. 
\end{itemize}

 But for quantum $SL(4, \Cb)$'s the classification  is unknown. 
 
 Here is another approach to cosemisimple Hopf algebras.

\begin{itemize}[leftmargin=*]
\medbreak\item    A  Hopf algebra $H$ is \emph{reductive} 
if  $\rep H$  is a semisimple category. \footnote{The $*$-Hopf algebras $H$
	such that any  finite-dimensional $H$-module carries an invariant inner product are reductive.
}
\end{itemize}

Thus $H^{\circ}$ and any of its Hopf subalgebras is cosemisimple. 
Reductive Hopf algebras which are pointed with abelian group, infinitesimal braiding of diagonal type 
and other suitable hypotheses were classified in \cite{ARS}.

\subsection{The lifting method II}\label{subsec:Method2}
The following approach was proposed in   \cite{A-Cuadra}.
Suppose that the coradical $H_0$ of $H$ is not a Hopf subalgebra of $H$. 
Then we look for a subcoalgebra $D$ of $H$ such that 
the $D$-wedge filtration $\left(\wedge^{n} D\right)_{n \in \N_0}$ is exhaustive and a Hopf algebra filtration.
Clearly $D$ should be a subalgebra and we already know that $D$ should contain $H_0$ \cite[5.2.8]{Montgomery-book},
thus the simplest choice is $D = H_{[0]} \coloneqq\ku\langle H_0\rangle$, which is called the Hopf coradical of $H$. 
We set $ H_{[n]}  \coloneqq \wedge^{n} D$ and  call this the standard filtration of $H$.
Again the corresponding graded coalgebra $\gr H$ is  a Hopf algebra that splits as the  bosonization
$\gr H \simeq \Bc \# H_{[0]}$
where $\Bc = \oplus_{n\ge 0} \Bc^n$ is a graded connected Hopf algebra in $\yd{H_{[0]}}$.
Thus we can face the analogues of  the steps \ref{item:levante0}, \dots, \ref{item:levante3} in this context.
Some features of this situation:
\begin{enumerate}[leftmargin=*,label=\textrm{(\alph*)}] 
\medbreak\item  \label{item:levanteII-0} There is no method to classify non-cosemisimple Hopf algebras generated by the coradical and having finite dimension or $\GK$. 
\medbreak\item \label{item:levanteII-1} 
The subalgebra of $\Bc$ generated by $V = \Bc^1$ 
is just  a pre-Nichols algebra of $V$. 
\end{enumerate}

Recall that a  finite-dimensional  algebra is basic 
when every irreducible representation has dimension 1; 
thus basic Hopf algebras are duals of pointed ones.  
Basic Hopf algebras are not necessarily generated by its coradical but among them 
there are many with this property.
In this context there is a fairly complete analysis. 

\begin{theorem}\label{thm:basic} \cite{A-Angiono-basic}
Let $L$ be a basic Hopf algebra generated by its coradical
such that $G = \Hom_{\text{alg}} (L, \ku)$ is abelian.
Let $\Htt$ be a  Hopf algebra  with $\Htt_{[0]} \simeq L$, so that
$\gr \Htt \simeq \Bc \# L$.
Then $\dim \Htt < \infty \iff \dim  \Bc < \infty \iff \Bc \simeq \toba(Z)$, where 
\begin{align*}
Z = \Gc(L(\lambda_1)) \oplus \dots \oplus \Gc(L(\lambda_t)) \in \yd{L}
\end{align*}	
is semisimple,  with $\lambda_1, \dots, \lambda_t \in \Irr \yd{\ku G}$, and 
$\dim \toba (V \oplus \lambda_1 \oplus \dots \oplus \lambda_t) < \infty$.
\end{theorem}

Here $V \in \yd{\ku G}$ is the infinitesimal braiding of $L^*$;
the $L(\lambda)$, $\lambda \in \Irr \yd{\ku G}$, are the simple objects in $\yd{\gr L^*}$
and $\Gc: \yd{\gr L^*} \to \yd{L}$
is a braided tensor equivalence whose definition requires \cite{Angiono-Garcia}. The proof, via the splitting technique in page \pageref{text:splitting-technique}, 
requires a large amount of the theory of pointed Hopf algebras with abelian group.

\subsection{The lifting method III}\label{subsec:Method-radical}
The lifting methods discussed in \S \ref{subsec:Method} and \S \ref{subsec:Method2} start from suitable coalgebra filtrations. Here we speculate about potential methods from algebra filtrations. 
That is,  we look at Hopf ideals $I$ such that the filtration $\Fc^I$  
is Hausdorff \eqref{eq:hausdorff}, i.e.,
\begin{align}\label{eq:hausdorff-I}
\bigcap_{n \in \N} I^n &= 0.
\end{align}
Then $\Fc^I$ is a Hopf algebra filtration and $\Gr_I H$ is a graded Hopf algebra 
that splits as the  bosonization
$\Gr_I H \simeq \Bc \# H/I$
where $\Bc = \oplus_{n\ge 0} \Bc^n$ is a graded connected Hopf algebra in the category $\yd{H/I}$.
By definition $\Bc$ is generated as an algebra by $V = \Bc^1$; i.e., $\Bc$ is a pre-Nichols algebra of $V$.
Thus we may consider the analogues of  the steps \ref{item:levante0}, \dots, \ref{item:levante3} in this context.
To turn this on, we have to single out a class of  Hopf ideals that satisfy \eqref{eq:hausdorff-I}. 
When $H$ is finite-dimensional,  the filtration $\Fc^{\Jac H }$ 
is  orthogonal to the coradical filtration of $H^*$. This suggests a first candidate, the Jacobson radical;
thus $H/ \Jac H $ is semiprimitive. We need to assume:
\begin{multicols}{2}
\begin{enumerate}[leftmargin=*,label=\textrm{(\roman*)}]
\medbreak\item \label{item:levanteIII-1} $\Jac H $ is a Hopf ideal,
\medbreak\item \label{item:levanteIII-2} \eqref{eq:hausdorff-I}, i.e., $ \bigcap_{n \in \N} (\Jac H )^n = 0$.
\end{enumerate}
\end{multicols}
Perhaps \ref{item:levanteIII-2} is redundant, see Question \ref{question:jacobson-hopf}. \phantomsection\label{text-question:jacobson-hopf}
Anyway, Question \ref{question:semiprimitive} appears naturally.
A second candidate might be the  ideal $\rsd H = \bigcap_{V \in \rep H} \ker \rho_V = \ker \iota$,
where $\iota: H \to (H^{\circ})^{\circ}$ is the natural map;
thus $\rsd  H $ is a Hopf ideal and $H/ \rsd  H$ is residually finite-dimensional.  
We need again  to assume the analogue of \ref{item:levanteIII-2}, cf. Question \ref{question:residue-hopf}. \phantomsection\label{text-question:residue-hopf}
In conclusion, this approach remains largely unexplored.

\section{Questions}\label{section:questions}
Throughout $H$ is a Hopf algebra. At the end of each Question we refer to the page where it is discussed, if so.

\begin{question}  \label{question:fflat-subalg} \cite{Montgomery-book}
If the antipode of $H$ is bijective,  is $H$  faithfully flat over any Hopf subalgebra with bijective antipode? 
P. \pageref{text-question:fflat}
\end{question} 

\begin{question}  \label{question:fflat-subalg-Noetherian} \cite{Skryabin-antipode-Noetherian}
If $H$ is Noetherian,  is it  residually finite-dimensional? 
P. \pageref{text-question:fflat}
\end{question}

\begin{question} \cite{Zhuang}
Does $\GK H$ belong to $\Z \cup \infty$?
\end{question}

The next two questions are classical open problems.

\begin{question} \label{question:enveloping-Noetherian}
If $U(\g)$ is Noetherian, is the Lie algebra $\g$ finite-dimensional? P. \pageref{text-question:enveloping-Noetherian}
\end{question}

\begin{question} \label{question:group-algebra-Noetherian}
If $\ku \varGamma$ is Noetherian,  is then $\varGamma$
 polycyclic-by-finite? P. \pageref{text-question:group-algebra-Noetherian}
\end{question}

\begin{question}\label{question:Noetherian-implies-affine} \cite{Wu-Zhang}
If $H$ is Noetherian,  is it necessarily affine?  P. \pageref{text-question:Noetherian-implies-affine}
\end{question}

\begin{question}\label{question:affine+finiteGK-Noetherian} \cite{Brown-Gilmartin-survey}
If $H$ is affine and $\GK H < \infty$, is $H$ necessarily  Noetherian? P. \pageref{text-question:affine+finiteGK-Noetherian}
\end{question}

\begin{question}  \label{question:subalg-Noetherian} 
If $H$ is Noetherian,  is  any Hopf subalgebra also Noetherian?  P. \pageref{text-question:subalg-Noetherian}
\end{question}

\begin{question}  \cite{Brown-survey} If $H$ is affine and PI, is it necessarily  Noetherian?
\end{question}

\begin{question}\label{question:cofrob-coss}
Classify the affine co-Frobenius, or cosemisimple, Hopf algebras that either have finite Gelfand-Kirillov dimension, or else are Noetherian. P. \pageref{text-question:cofrob-coss}
\end{question}

\begin{conjecture}\label{conjecture:antipode} \cite{Skryabin-antipode} 
 The antipode of a Noetherian Hopf algebra is bijective. P. \pageref{subsec:Antipode}
\end{conjecture}

\begin{question} \cite{Brown-survey,Brown-Zhang}  If $H$ is affine Noetherian and PI, 
does $\Ss$ have finite order? Is some (even) power of $\Ss$ an inner automorphism of $H$?
\end{question}

\begin{conjecture}\label{conjecture:0-divisor} (Kaplansky)\footnote{This question appeared earlier in \cite{Higman} where it is proved for locally indicable groups.}
If $\varGamma$ is torsion-free, then $\ku \varGamma$ is a domain. P. \pageref{text:conjecture:0-divisor}
\end{conjecture}

\begin{question}\label{question:classification:domain-prime-semiprime}
Classify affine  Hopf algebras with finite $\GK$ that are domains, prime, or  semiprime.
P. \pageref{text-question:classification:domain-prime-semiprime}
\end{question}

\begin{question}\label{question:domain-semiprimitive} (Skryabin, personal communication).
If $H$ is a Noetherian domain, is it forcely semiprimitive, i.e., $\Jac H = 0$? P. \pageref{text-question:domain-semiprimitive}
\end{question}

\begin{question} \cite{Brown-survey}  If $H$ is  a Noetherian PI domain, is it necessarily  
	module-finite over its centre? What if it is also affine?\footnote{The original version  with the hypothesis \emph{prime} instead of domain was disproved in \cite{Gelaki-Letzter}.}
\end{question}

\begin{question} \cite{Wang-ZZ-JAlg}  If $H$ is affine, Noetherian and a domain, does
$\GK H \leq 2$ imply that it is generated by grouplikes and skew primitives?
\end{question}

\begin{question}\label{question:nichols-locally-finite} \cite{A-Angiono-Heckenberger-infinite}
Let $(V, c)$ be a braided vector space with $\GK \toba(V) < \infty$.
Is $(V, c)$ a locally finite braided vector space? P. \pageref{text-question:nichols-locally-finite}
\end{question}

\begin{question}\label{question:H-affine-grH-affine}
If $H$ is affine, $H_0$ is a Hopf subalgebra and $\GK H < \infty$, is $\gr H$ affine?
P. \pageref{text-question:H-affine-grH-affine}
\end{question}

\begin{question}\label{question:postnichols-finite-gkd-gral}
Determine  $\postfh{H}(V)$ for $V \in \yd{H}$ with $\GK \toba(V) < \infty$. P. \pageref{text-question:postnichols-finite-gkd-gral}
\end{question}

\begin{question}\label{question:postnichols-finite-gkd-diagonal}
Determine $\pref(V)$ when $V$ is of diagonal type with braiding matrix  of Cartan type $A_{\theta}$ or $D_{\theta}$ with $q=-1$. P. \pageref{text-question:postnichols-finite-gkd-diagonal}
\end{question}

\begin{question}\label{question:nichols-finite-gkd} Let $\varGamma$ be a  nilpotent-by-finite group.
Determine all post-Nichols algebras  $\Rc$ of  $V \in \yd{\ku \varGamma}$ such that $\GK \Rc < \infty$.
P. \pageref{text-question:nichols-finite-gkd}
\end{question}

\begin{question}\label{question:racks-finite-gkd}
Determine the pairs  $(X, \bqg) $,  $X$  a \emph{finite} rack   and $\bqg: X\times X \to GL(n, \ku)$
a 2-cocycle, such that $\GK \toba(X, \bqg) < \infty$. P. \pageref{text-question:nichols-finite-gkd}
\end{question}

\begin{conjecture}\label{conjecture:typeC-X-finite-GK} \cite{A-nilpotent}
Let $\Orb$ be a finite rack  of type C.
Then $\GK \toba(\Orb, \bqg) = \infty$ for every faithful 2-cocycle $\bqg$. P. \pageref{text-conjecture:typeC-X-finite-GK}
\end{conjecture}

\begin{question}\label{question:nichols-cocycle-1}  
Let $X$ be an  indecomposable rack with $1 < \vert X \vert < \infty$
and let $\mathbf{1}: X \times X \to \kut$ be the constant cocycle  1;
is $\GK\toba(X, \mathbf{1}) < \infty$? P. \pageref{text-question:nichols-cocycle-1}
\end{question}

\begin{question}\label{question:limsup} Let $\Bc$ be a pre-Nichols algebra of a 
finite-dimensional braided vector space  $(V,c)$. Assume  that $\GK \Bc < \infty$.
 Is $\GK \Bc =  \lim_{n \to \infty} \log_n \dim \Bc^n$?
Is the Hilbert series of $\Bc $ a rational function?
\end{question}

\begin{question}\label{question:braided-tensor-product}
Is $\GK \toba(V \underline \otimes W) = \GK \toba(V) + \GK \toba(W)$ provided that
$V, W\in \yd{H}$ satisfy  $c_{V, W}c_{W, V} = \id_{W \otimes V}$?
P. \pageref{text-question:braided-tensor-product}
\end{question}

\begin{question}\label{question:CQG}
Classify the affine compact Hopf algebras that either are Noetherian, or else 
have finite $\GK$.
\end{question}

\begin{question}\label{question:CQG-Lie} \cite{Woronowicz-newquantum}
Classify the CQG, or more generally the cosemisimple Hopf algebras, of Lie type.  P. \pageref{def:CQG}
\end{question}

\begin{question}\label{question:CQG-Lie-exhaustion} Is a simple (i.e., without normal Hopf subalgebras) cosemisimple Hopf algebra $H$ with  $0 < \GK H < \infty$ necessarily of Lie type?
P. \pageref{text-question:CQG-Lie-exhaustion}
\end{question}

\begin{question}\label{question:jacobson-hopf}
If  $\Jac H $ is a Hopf ideal, when does it satisfy $ \bigcap_{n \in \N} \Jac H ^n = 0$? P. \pageref{text-question:jacobson-hopf}
\end{question}

\begin{question}\label{question:semiprimitive}
For $H$ semiprimitive with $\GK H < \infty$, 
classify all  pre-Nichols algebras $\Bc$ in $\yd{H}$ with $\GK \Bc < \infty$, and
all their liftings, i.e.,  all Hopf algebras $L$ such that $\Gr L \cong \Bc\# H$.
P. \pageref{text-question:residue-hopf}
\end{question}

\begin{question}\label{question:residue-hopf}
When does the ideal $\rsd  H = \bigcap_{V \in \rep H} \ker \rho_V$  satisfy \eqref{eq:hausdorff-I}?
P. \pageref{text-question:residue-hopf}
\end{question}

\begin{ack} I thank the Scientific Committee of ICRA 2022 for the invitation to give a conference 
and to present this contribution. I am   grateful to
Iv\'an Angiono, Julien Bichon, Ken Brown, Juan Cuadra, 
Fran\c{c}ois Du\-mas, Pavel Etingof, Ken Goodearl, 
Istv\'an Heckenberger,  Sonia Natale, Héctor Martín Peña Pollastri,
Guillermo Sanmarco, Peter Schauenburg, Hans-J\"urgen Schneider, Serge Skryabin,  Blas Torrecillas,
Milen Yakimov, Hiroyuki Ya\-mane  and James Zhang
for  fruitful discussions, exchan\-ges and collaborations.

\end{ack}

\end{document}